\documentclass{amsart}

\usepackage{graphicx}
\usepackage{float}
\usepackage{amsmath}
\usepackage{amsfonts}
\usepackage{graphicx}
\usepackage{subfigure}
\usepackage{cite}
\usepackage[colorlinks,linkcolor=blue,urlcolor=blue,citecolor=blue]{hyperref}
\usepackage{todonotes}
\usepackage{siunitx}

\renewcommand{\vec}[1]{{\mathbf #1}}
\newcommand{\R}{{\mathbb R}}

\newtheorem{remark}{Remark}
\newtheorem{proposition}{Proposition}

\begin{document}

\title{Estimating and using information in inverse problems}

\author{Wolfgang Bangerth}
\address{Department of Mathematics, Department of Geosciences, Colorado State University,
  Fort Collins, CO 80523-1874, USA}
\email{bangerth@colostate.edu}
  
\author{Chris R. Johnson}
\address{Scientific Computing and Imaging Institute, University of Utah,
Salt Lake City, UT  84112, USA}
\email{crj@sci.utah.edu}

\author{Dennis K. Njeru}
\address{Scientific Computing and Imaging Institute, University of Utah,\\
Salt Lake City, UT  84112, USA}
\email{dnjeru@sci.utah.edu}

\author{Bart van Bloemen Waanders}
\address{Sandia National Laboratories, MS 0370, PO Box 5800,
  Albuquerque, NM 87185, USA}
\email{bartv@sandia.gov}

\renewcommand{\shortauthors}{W. Bangerth, C. R. Johnson, D. K. Njeru,
  B. van Bloemen Waanders}

\begin{abstract}
  In inverse problems, one attempts to infer spatially variable
  functions from indirect measurements of a system.  To practitioners
  of inverse problems, the concept of ``information'' is familiar when
  discussing key questions such as which parts of the function can be
  inferred accurately and which cannot. For example, it is generally
  understood that we can identify system parameters accurately only
  close to detectors, or along ray paths between sources and
  detectors, because we have ``the most information'' for these
  places.
  
  Although referenced in many publications, the ``information'' that
  is invoked in such contexts is not a well understood and clearly defined
  quantity. Herein, we present a definition of \textit{information
    density} that is based on the variance of coefficients as derived
  from a Bayesian reformulation of the inverse problem. We then
  discuss three areas in which this information density can be useful
  in practical algorithms for the solution of inverse problems, and
  illustrate the usefulness in one of these areas -- how to choose the
  discretization mesh for the function to be reconstructed -- using
  numerical experiments.
\end{abstract}


\maketitle

\section{Introduction}

Inverse problems -- i.e., determining distributed internal parameters of a
system from measurements of its state -- are frequently
ill-posed. Mathematically, this ill-posedness is often described as the lack
of a continuous mapping from the space of measurements to the corresponding
parameters reconstructed from a measurement. A consequence of ill-posedness is that a small
measurement error can result in a significantly different reconstructed
parameter unless the problem is \textit{regularized} in some way.

Concretely, let us consider that we want to identify a spatially varying
parameter $q=q(\mathbf x)$, for instance the density and elastic moduli of the earth
in seismology or the absorption and scattering properties of the human body
in biomedical imaging. This identification requires using measurements $z$ of some part of the state of
the system under interrogation, e.g. the time-dependent displacement at a
seismometer station, or the light intensity at the surface of the body as
recorded by the pixels of a camera. If $z$ is corrupted by noise of level
$\varepsilon$ we will get two reconstructions ($q_1=q_1(\mathbf x),q_2=q_2(\mathbf x)$) for two measurements
($z_1,z_2$) that differ only by the realization of the measurement
noise. Ideally, we would be able to show that
\begin{equation}
  \label{eq:stability}
  \| q_1-q_2 \| \le C \| z_1-z_2 \|
\end{equation}
for some appropriate choice of norms and a constant $C$ of
moderate size. The problem is ``ill-posed'' if such
an estimate does not exist. Many inverse problems fall into this category of ill-posedness.

On the other hand, a pragmatic view of inverse problems is that the
ill-posedness of the problem is simply the result of a \textit{lack of
  information}. Some inverse problems can achieve well-posedness by obtaining
different kinds of measurements from the system under consideration, but even
if that is not the case, it would already be advantageous to simply
\textit{reduce the ill-posedness}. In either case, we may ask whether it
is possible to derive estimates of the kind
\begin{equation}
  \label{eq:information-bound}
  \| J(\mathbf x) (q_1(\mathbf x)-q_2(\mathbf x)) \| \le C \| z_1-z_2 \|,
\end{equation}
again with an appropriate choice of norms. We will call $J(\mathbf x)\ge 0$,
or a related quantity such as its square root, the \textit{information
  density}. Equation (\ref{eq:information-bound}) is motivated by the observation that at
places where $J(\mathbf x)$ is large, we can accurately determine
the value of the coefficient $q(\mathbf x)$ we are looking for (i.e.,
$|q_1(\mathbf x)-q_2(\mathbf x)|$ must be small to satisfy the inequality). Conversely,
the places where $J(\mathbf x)$ is small coincide with those locations
where we have little control over the coefficient, and even small amounts of
noise in $z$ may lead to large variations $|q_1(\mathbf x)-q_2(\mathbf
x)|$. In the extreme case when the problem is truly ill-posed, 
$J(\mathbf x)$ would not be bounded away from zero and consequently $\| J(\mathbf x)
\varphi(\mathbf x) \|$ would not be a norm of $\varphi$.%
\footnote{In many inverse problems, for
  example in imaging, the ill-posedness manifests not by there being
  \textit{locations} $\mathbf x$ at which $q(\mathbf x)$ is not
  identifiable. Rather, it is the \textit{high-frequency (Fourier)
    content} of $q$ that is often not identifiable without
  regularization. In this case, an equation
  like \eqref{eq:information-bound} with $\mathbf x$ replaced by the
  wave number $\mathbf k$ can be considered. Regardless of this obvious difference, we move forward with the derivation as stated.}

It is unlikely that for practical problems we can find meaningful
expressions for $J$ that give rise to provable estimates of the form
\eqref{eq:information-bound}. This is because for many inverse
problems, what can and cannot be recovered stably is often not about
where in space we are, but about which \textit{modes in feature space}
(for example low- versus high-frequency components of a function
$q(\mathbf x)$) are identifiable.
In our discussions below, we therefore consider estimates such as
\eqref{eq:information-bound} \textit{aspirational}: We will instead seek
statements such as
\begin{equation}
  \label{eq:information-bound-approximation}
  \| j(\mathbf x) (q_1(\mathbf x)-q_2(\mathbf x)) \| \simeq C \| z_1-z_2 \|,
\end{equation}
where $j(\mathbf x)$ takes on the role of the information density, and
where $\simeq$ expresses a relationship of the form ``behaves conceptually
like, but possibly only when spatial discretization is used''. Although proving 
that any choice of $j$ in \eqref{eq:information-bound-approximation} implies
an estimate of the form \eqref{eq:information-bound} is likely impossible, the conceptual approach of seeking a
function $j(\mathbf x)$ that expresses the idea of an information
density of how much we know about $q$ at different points in space will turn out to be useful in practice -- as we will
demonstrate in Sections~\ref{sec:examples-of-use} and \ref{sec:results}.

\paragraph{References to information in inverse problems in the research literature}
The notion of information density is not new, in particular in applications where $q(\mathbf x)$ is replaced by finite-dimensional parameter vectors $(q_i)_{i=1}^N$. Indeed,  similar notions can be found in
many areas of inverse and parameter estimation problems in various
forms, and among practitioners of inverse problems, there is a degree of ``knowledge'' that \textit{information} is a key concept. At the same time, practitioners do not appear to have a clear understanding of what information actually means, and
uses of this concept in the literature appear to be vaguely defined and disconnected. In practice, references to the term ``information'' in the literature on inverse problem are almost always qualitative, rather than giving the term a quantitative definition.

Yet, many publications touch on the \textit{concept} of information in inverse problems. The most obvious application of information concepts to inverse problems is in optimal experimental design where the goal of the design of schemes is to measure data about the system to minimize the uncertainty (that is, to maximize the information) in the parameters we wish to recover \cite{AD92,BA11,EM98}. This relation to uncertainty is most clearly articulated in the Bayesian setting of
optimal experimental design where the
\textit{information gain} of the posterior probability distribution over the prior is maximized.  However, information, generally defined in less concise terms, is also a topic discussed in other contexts. For example, considering concrete applications,
\cite{RO05} presents Fisher information for a single-particle system and proposes
a new uncertainty relationship based on Fisher information. Similarly, \cite{Zh07,BE00a} discuss
the use of a resolution matrix in seismic tomography (see also \cite{RHDW11}); related concepts of resolution, resolution length scales, event kernels, sensitivity kernels, Fr{\'e}chet kernels, or point spread functions also appear in both seismic imaging and a number of other fields, see for example \cite{Kuc09,DXGKL11,VKK00,Liu_2012,Montelli_2006,Rawlinson_2014,GLAD-M15,GLAD-M25,Zhao2005,Dahlen2000,Hung2000}.
In many other cases, the literature references the Fisher information matrix that, together with the Cram{\'e}r-Rao bound, quantifies how accurately we know what the inverse problem seeks to identify \cite{Kay_93}; examples include \cite{BE10}, which uses this approach
for estimating diffusion in a single particle tracking process; \cite{FA07},
which compares Fisher matrices to the Hessian calculation in boundary value
inversion problem using the heat equation; and \cite{NO09} which presents a preconditioning
and regularization scheme based on Fisher information. 

These publications generally state that there is a connection between the number, kind, and
accuracy of measurements on the one hand, and the uncertainty in the recovered
parameters of the inverse problem on the other. But, none of the studies mentioned go on to specifically \textit{identify the role of information in the spatially variable ability to recover parameters in inverse problems} in a systematic way. Let us illustrate this with two publications that use data from seismometer recordings of earthquakes to produce reconstructions of the wave speeds within the Earth. Specifically, \cite{GLAD-M15} presents the GLAD-M15 earth model that is based on 253 earthquake events and recorded at seismometer stations around the world; \cite{GLAD-M25} is an improvement over the previous paper and presents the GLAD-M25 model that is based on 1040 events. Because neither earthquakes nor seismometers are uniformly distributed across the earth surface, there are parts of the earth where we have a fairly good understanding of how fast seismic waves travel (namely those parts that are intersected by many ``ray paths'' from earthquake sources to seismometer stations) whereas there are other parts we do not really know very much about. In other words, we have ``much information'' about some parts and ``little information'' about others; such statements are certainly common in many inverse problems, but information is not typically a quantitative term in these statements. Indeed, the second paper mentions the term ``information'' eight times on 21 pages, always in qualitative contexts such as ``The ultimate goal is to use every single piece of information in seismograms'' and ``Figs 2–4 illustrate that much information in seismograms is being assimilated''. Yet, the paper does not present a quantitative definition of ``information''. The authors do, however, attempt to quantify the accuracy of their reconstructions: Both \cite{GLAD-M15} and \cite{GLAD-M25} compute ``point-spread functions'' that assess the size of the smallest features in the earth the method can resolve given the available seismic observations. As the papers illustrate well, the size of the smallest resolvable features varies substantially between different parts of the earth, in accordance with our statements about having more or less information depending on location. Furthermore, the second paper illustrates the improvement that results from using more seismic events by showing that the smallest resolvable features have become smaller in GLAD-M25 compared to GLAD-M15 \cite[Section 6.1 and Fig. 7]{GLAD-M25}. However, the evaluations performed there are entirely qualitative: The paper states that ``We observe that the 100 km diameter $\beta_V$ anomaly---which is smaller than what was used in GLAD-M15---at 300 km depth is well recovered'', and ``In each case the Gaussian anomaly is well resolved'' (end of Section 6.1). It is these kinds of comparisons that we would like to put onto a quantitative basis, and the method we will discuss herein allows us to do exactly this by providing an ``information density'' that (i)~would be able to quantify how much more information is available in some parts of the world than in others, and (ii)~that would allow for a quantitative evaluation of the improvement of the second over the first model.

One could similarly analyze papers from many other disciplines that use inverse problems. They may be using different words, but a common feature of the many definitions of resolution, adjoints, sensitivity, and identifiability that can be found in the literature, is that most of these notions originate in, and were developed for \textit{deterministic} inverse problems. On the other hand, ``information'' is probably best understood as a statistical concept, and a useful definition will therefore be rooted in statistical reformulations of the inverse problems. We will utilize the connection between the Fisher information matrix and the variance of reconstructed parameters, via the Cram{\'e}r-Rao bound, to derive information densities in Section~\ref{sec:finite-dimensional-information} below. If the parameter-to-measurement map is linear, our definition of ``information'' is related to the sensitivity of this map to changes in parameters; in the nonlinear case, it is related to some average of this sensitivity.

The differences in the concepts mentioned above, and the lack of a common language to describe them, presents the motivation for this work. Indeed, while many of the referenced works above at least verbally express the notion that in many inverse problems, we ``know more'' about the parameters in some parts of the domain than in others, we have not found published works that try to provide quantitative measures of this concept. The current work tries to address this deficit; in concrete terms, our goals are as follows:
\begin{itemize}
    \item To introduce the notion of an \textit{information
  density} based on a statistical interpretation of the inverse problem, the Fisher information matrix, and an application of the Cram{\'e}r-Rao bound.
  \item To define the information density in a way that respects intuitive notions. In particular, it should satisfy that adding more measurements can only increase the information, not decrease it; that information is inverse proportional to measurement uncertainties; and that with Gaussian noise and linear models, measuring the same quantity twice increases the information by $\sqrt{2}$.
  \item To outline a number of applications for which we believe information density can be usefully employed.
  \item To practically evaluate our concepts in a concrete application, namely the choice of mesh on which to discretize an inverse source-identification problem.
\end{itemize} 

Herein, let us first provide some perspectives in Section~\ref{sec:perspective} on how one formulates inverse problems, and how the different philosophical approaches inform our approach. We then address the goals mentioned in the previous paragraph by first considering a finite-dimensional, linear model problem in Section~\ref{sec:finite-dimensional} that we use to provide a conceptual overview of what we are trying to achieve, followed by the extension of this model problem to the infinite-dimensional case in Section~\ref{sec:infinite-dimensional}. Having so set the stage, in Section~\ref{sec:examples-of-use}, we provide ``vignettes'' for three ways in which we believe information densities can be used in practice. Section~\ref{sec:results} then explores one of these -- the choice of mesh for discretizing an infinite-dimensional inverse problem -- in detail and with numerical and quantitative results. We conclude in Section~\ref{sec:conclusions}. Two appendices discuss the extension of our work to nonlinear problems (Appendix~\ref{sec:nonlinear}) and explain the derivation of the mesh refinement criteria that we compare against the method we propose in Section~\ref{sec:results} (Appendix~\ref{sec:amr-ee}).

\section{Perspectives on inverse problems}
\label{sec:perspective}

Parameter estimation problems -- of which inverse problems are a particular kind -- seek information about parameters in the model using measurements (``observations'') of the system's state. Depending on what kind of, and how much, information one seeks, parameter estimation problems can be formulated in different ways. For instance,   \textit{deterministic parameter estimation problems} only seek very limited information about parameters (namely, a ``most likely'' value) and in return can be solved for large numbers of inversion parameters with complicated physics.  On the other hand, \textit{Bayesian parameter estimation problems} are formulated to achieve complete statistical characterization but are much more expensive to solve and therefore limited to relatively small parameter spaces and simpler physics. Deterministic and Bayesian inverse problems are therefore often considered the extremes on a spectrum of formulations.



Historically, parameter estimation problems were usually formulated as seeking that set of parameters $q^\ast$ for which model predictions fit measurements \textit{best}. This approach is often called the \textit{deterministic} approach to parameter estimation. Oftentimes, one measures fit of predictions to observations using the $l_2$ norm, in which case the problem becomes a least-squares problem. A typical example every undergraduate learns about is to find the parameters $q=\{a,b\}$ in a linear model $y=ax+b$ that best fit a set of data points $\{x_i,y_i\}_{i=1}^N$. If the parameter we seek is a distributed function $q(\mathbf x)$ that appears in a partial differential equation, then the inverse problem is often formulated as minimizing the difference between model predictions and measured data, subject to the PDE as a constraint that connects parameters and state variables from which one can then extract predicted measurements. These so-called PDE-constrained optimization problems can be solved relatively efficiently for high-dimensional inversion spaces, requiring perhaps a few hundred solutions of the model even if the parameter-to-measurement map is nonlinear \cite{BBGMSW12,Noemi2014,Ban08ip,Akcelik2006}.
As a consequence, many practically relevant inverse problems can be solved using this approach, including among many others some biomedical imaging modalities \cite{BJ07ip} and seismic tomography in which one seeks to determine the Earth's three-dimensional geologic make-up based on measured responses of Earth to earthquakes or artificial signals \cite{GLAD-M25}.

The other end of the spectrum of parameter estimation formulations is occupied by an approach commonly referred to as \textit{Bayesian}, a perspective on inverse problems that for many in our community was first popularized by Albert Tarantola's 1987 book \cite{Tar87} and its 2004 re-issue by SIAM \cite{Tar04}. In this approach, one does not seek just the \textit{best} parameter $q^\ast$, but in fact a probability distribution $p(q|z)$ that, in intuitive terms, describes how likely a parameter $q$ is given what we have measured (measurements being denotes by $z$). Because one can only obtain $p(q|z)$ \textit{after} measuring, it is called the \textit{posterior} probability distribution. It takes into account that all measurements we can produce are always uncertain, and that perhaps also the model is uncertain or inaccurate. Clearly, having a whole probability distribution is much more information than just a single, best estimate of a parameter. Among the benefits of the Bayesian approach is that we can quantify the uncertainty in the estimated parameters, for example by computing the variance of parameters under $p(q|z)$, along with correlations between parameters via the covariances of $p(q|z)$. Typically, the formulations of deterministic and Bayesian estimation problems also imply that the best fit (deterministic) estimate equals the point where $p(q|z)$ is maximal -- i.e., $q^\ast$ is the \textit{maximal a posteriori} (MAP) point of $p(q|z)$. 

The Bayesian perspective provides a much richer description of parameter estimation problems, but it comes at a cost: $p(q|z)$ is in general a probability distribution that has no closed form formula, and that can only be explored by statistical sampling -- usually using Monte Carlo Markov Chain (MCMC) methods. Each sample typically requires the solution of a forward model, and in practical inverse problems one therefore needs thousands, millions, or even more PDE solves. (For an extreme case, see \cite{AB23} in which highly accurate statistics are computed by solving the PDE $2\cdot 10^{11}$ times.) In other words, to compute good statistics for Bayesian inverse problems is \textit{orders of magnitude more expensive} than solving deterministic ones. These cost considerations have so far prevented fully Bayesian formulations from being used widely for practical, large-scale inverse problems.

\begin{figure}
    \centering
    \includegraphics[height=0.44\textwidth]{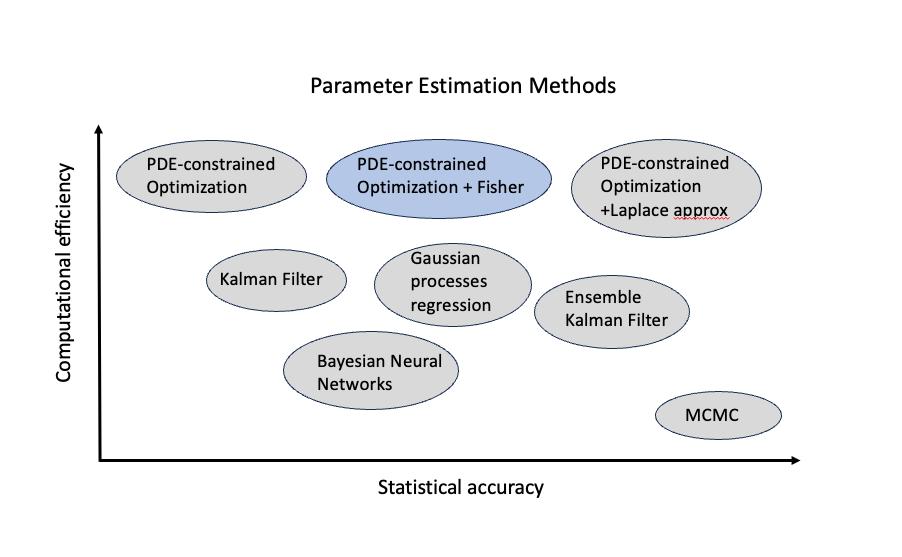}
    \caption{\it Notional positioning of parameter estimation methods. Our approach -- enriching deterministic parameter estimation problems solved as a PDE-constrained problem with information obtained from the Fisher information matrix -- is shown in blue.}
    \label{fig:efficiency_accuracy}
  \end{figure}

As a consequence, many methods have been developed over the years that use insights from deterministic problems to make the Bayesian problem more tractable, or to \textit{approximate} the latter. The result are methods that lie somewhere on the spectrum between the deterministic and Bayesian end points; for example, one could try to approximate $p(q|z)$ by a Gaussian in return for substantial computational savings. Figure \ref{fig:efficiency_accuracy} provides a notional characterization of how different well-known methods could be positioned on a spectrum that takes into account statistical accuracy and computational efficiency.  The extremes of the spectrum are occupied by PDE-constrained optimization (used for deterministic optimization problems: computationally efficient but statistically not very insightful) and MCMC solvers for Bayesian inference (statistically accurate but computationally very expensive).  In between, a range of hybrid methods can be positioned.   The comparison is of course not straightforward but \textit{conceptual}; for example, some information can be computed online or offline, and some methods require sequential computations whereas others can do things simultaneously in parallel. 
Nevertheless, the characterization provides a general perspective of utilizing information in inverse problems.

The methods in the middle of the spectrum mentioned above and in the figure use deterministic methods to make the solution of statistical formulations more efficient.
Yet, relatively little has been done to bring information from the (expensive) Bayesian perspective to selectively enrich the deterministic solution, at a cost that is comparable to solving the deterministic problem. Indeed, we are not aware of any literature that would have done so in the context of inverse problems, i.e., for determining spatially variable coefficients in partial differential equations despite the fact that our literature review in the introduction has found much interest in the topic. This paper then should be seen in this context: While we seek a best-fit (deterministic) estimate of a spatially variable parameter in a partial differential equation, we would like to enrich the information we have about these parameters by a spatially variable measure of \textit{information} that can help us with tasks such as choosing the mesh on which to discretize the parameter. In order to define this \textit{information density}, we will resort to the Bayesian perspective (but without actually solving a Bayesian inverse problem), using the fact that statistical analysis has provided us with ways of estimating standard deviations of parameter estimates at far lower cost than computing the probability distribution $p(q|z)$. In fact, our method only requires a number of model solutions proportional to the larger of the number of parameters and the number of observations. In the next section will then outline how this can be done for a finite-dimensional problem, and Section~\ref{sec:infinite-dimensional} illustrates the application of these ideas to an actual inverse problem in which the coefficient $q=q(\mathbf x)$ is a spatially variable function.

\section{A finite-dimensional, linear model problem}
\label{sec:finite-dimensional}
In order to explain our ideas, let us first consider a linear,
finite-dimensional problem. (In practice, of course, many inverse
problems are nonlinear; as we will show in
Appendix~\ref{sec:nonlinear}, it turns out that nearly everything we
say here will carry over to the nonlinear case under common
assumptions, namely, if the model is not too nonlinear, or if the noise level in measurements is sufficiently small.)
Specifically, let us consider a problem in which a state vector $\mathbf u$ is related to source terms via the relationship
\begin{equation}
    \label{eq:finite-dimensional-state-eq}
    A \mathbf u = \sum_k q_k \mathbf s_k, 
\end{equation}
where $\mathbf s_k$ are possible source vectors and $q_k$ are their relative strengths. We assume that the system matrix $A$ is invertible, although it may be ill-conditioned. In the inverse problem, we are then interested in recovering unknown source strengths $q_k$ through
a number of (linear, noisy) measurements
\begin{equation}
    \label{eq:finite-dimensional-measurements}
    z_\ell = \mathbf m_\ell^T \mathbf u + \varepsilon_\ell, 
\end{equation}
where the dot product $\mathbf m_\ell \cdot$ corresponds to the $\ell$th measurement operator, and
$\varepsilon_\ell$ is measurement noise. We assume that we have a guess $\sigma_\ell$ for the magnitude of the noise $\varepsilon_\ell$.

For convenience, let us collect the quantities $q_k,\mathbf s_k, z_\ell, \mathbf m_\ell$ into vectors and matrices $\mathbf q, S, \mathbf z, M$, where the $\mathbf s_k$ form the columns of $S$ and the $\mathbf m_\ell$ the rows of $M$. Then, we can state the source strength recovery problem we will consider here as
\begin{align}
    \label{eq:finite-dimensional-minimization}
    \begin{split}
    \min_{\mathbf q,\mathbf u} \; &
    \frac 12 \left\| M \mathbf u - \mathbf z \right\|^2_{\Sigma^{-2}} 
    + \frac{\beta}{2} \left\| R \mathbf q \right\|^2,
    \\
    \text{such that} \; &
    A \mathbf u = S \mathbf q.
    \end{split}
\end{align}
Here, we have used the weighted norm $\|M\mathbf u-\mathbf z\|_{\Sigma^{-2}}^2 = \sum_\ell \frac{1}{\sigma_\ell^2} |\mathbf m_\ell^T\mathbf u - z_\ell|^2$, where the diagonal matrix $\Sigma_{\ell\ell}=\sigma_\ell$ weighs measurements according to the assumed \textit{certainty} $\frac{1}{\sigma_\ell}$ we have of the $\ell$th measurement.
We have also added a Tikhonov-type regularization term where $\beta$ is the regularization parameter and $R$ a matrix that amplifies the undesirable components of $\mathbf q$.

By eliminating the state variable using the state equation \eqref{eq:finite-dimensional-state-eq}, we can re-state this problem as an unconstrained, quadratic minimization problem:
\begin{align}
    \label{eq:finite-dimensional-minimization-reduced}
    \begin{split}
    \min_{\mathbf q} \; {\mathcal{J}}_\text{red}(\mathbf q) :=&
    \frac 12 \left\| M A^{-1} S \mathbf q - \mathbf z \right\|^2_{\Sigma^{-2}} 
    + \frac{\beta}{2} \left\| R \mathbf q \right\|^2.
    \end{split}
\end{align}
It is then not difficult to show that the minimizer $\mathbf q$ of this problem satisfies
\begin{align}
    \label{eq:finite-dimensional-solution}
    \underbrace{
    \left(S^T A^{-T} M^T \Sigma^{-2} M A^{-1} S + \beta R^T R \right)}_Q \mathbf q = S^T A^{-T} M^T \Sigma^{-2} \mathbf z.
\end{align}

If we consider the noise to be random, we can ask how the solution $\mathbf q$ depends on concrete measurements. 
Specifically, if we have two, presumably nearby, measurements $\mathbf z_1, \mathbf z_2$, then the following relationship holds for the corresponding solutions $\mathbf q_1, \mathbf q_2$:
\begin{align}
    \label{eq:finite-dimensional-solution-difference}
    Q (\mathbf q_1 - \mathbf q_2) = S^T A^{-T} M^T \Sigma^{-2}(\mathbf z_1- \mathbf z_2).
\end{align}

For the following discussions, it is important to point out that the
structure of $Q$ guarantees that the matrix is symmetric, positive, and semidefinite; that is, all of its eigenvalues are non-negative. We will assume that the user has chosen $\beta$ and $R$ in such a way that $Q$ is positive definite, although some of the eigenvalues may be small.

\subsection{Defining an ``information content'' for components of the solution vector}
\label{sec:def-based-on-covariance}

We investigate herein if the relationship
\eqref{eq:finite-dimensional-solution-difference} between $\mathbf q_1
- \mathbf q_2$ and $\mathbf z_1 - \mathbf z_2$ allows us to define
stability bounds such as those outlined in
\eqref{eq:information-bound} or
\eqref{eq:information-bound-approximation} above.  In the context of this finite-dimensional situation, such a bound would have the form
\begin{align}
    \label{eq:finite-dimensional-bound}
    \left\|
      \mathbf j \odot (\mathbf q_1 - \mathbf q_2)
    \right\| 
    \simeq
    C
    \left\| \mathbf z_1- \mathbf z_2 \right\|,
\end{align}
with a constant $C$ of possibly unknown size, and
where $\odot$ indicates the Hadamard product that scales each entry of the vector $\mathbf q_1 - \mathbf q_2$ by the corresponding entry of the ``information vector'' $\mathbf j$. (Alternatively, we can interpret $\mathbf j \odot (\mathbf q_1 - \mathbf q_2)$ as $\text{diag}(\mathbf j) (\mathbf q_1 - \mathbf q_2)$ where $\text{diag}(\mathbf j)$ is a diagonal matrix with diagonal entries $j_k$.) 

A \textit{meaningful} statement%
\footnote{Equation
  \eqref{eq:finite-dimensional-solution} implies that
  $\|\mathbf q_1 - \mathbf q_2\| \le C \|\mathbf z_1 - \mathbf
  z_2\|$ with $C=\|Q^{-1} S^T A^{-T}M^T \Sigma^{-2}\|$, which
  corresponds to choosing $\mathbf j$ in
  \eqref{eq:finite-dimensional-bound} as a vector of ones. At the same
time, this estimate reflects no specifics of the problem and we do not
consider this choice useful because it does not help us identify \textit{which
components of $\mathbf q$ can be identified accurately, and which cannot}.}
such as \eqref{eq:finite-dimensional-bound} will not always follow
from \eqref{eq:finite-dimensional-solution} unless either the action
of the matrix $Q=S^TA^{-T}M^T \Sigma^{-2} M A^{-1} S + \beta R^TR$ can
somehow be approximated from below by a diagonal matrix, or $Q^{-1}$
be approximated from above by a diagonal matrix. Indeed, we could choose $\mathbf j$ to be a vector whose elements are all equal to $j_k := \lambda_\text{min}(Q)=[\lambda_{\max}(Q^{-1})]^{-1}$. If, in addition, $C:=\|S^T A^{-T} M^T\Sigma^{-2}\|$, 
then \eqref{eq:finite-dimensional-bound} holds true. This approach
works as long as the regularization is chosen so that all eigenvalues
of $Q$ are reasonably large, i.e., that the problem is well-posed; in
practice, however, this choice may over-regularize the problem. Alternatively, we could try to obtain a bound by formally solving 
\eqref{eq:finite-dimensional-solution-difference} to obtain
$    (\mathbf q_1 - \mathbf q_2) = Q^{-1} S^T A^{-T} M^T \Sigma^{-2}(\mathbf z_1- \mathbf z_2)
$, and then trying to find a vector $\mathbf j$ obtained from a lower bound of
\begin{align*}
  Q^{-1} S^T A^{-T} M^T \Sigma^{-2}
  =
  \left[S^T A^{-T} M^T \Sigma^{-2} M A^{-1} S + \beta R^T R\right]^{-1} \left[S^T A^{-T} M^T \Sigma^{-2}\right].
\end{align*}

In the end, we have found that choices such as $j_k = \lambda_\text{min}(Q)$ obtained by purely algebraic manipulation tend not to be
interesting since they do not give us any insight into \textit{which
  elements of $\mathbf q$ can be accurately estimated and which
  cannot}. Furthermore, if the problem is indeed ill-posed, as the
regularization parameter $\beta$ is reduced, \textit{all} elements
$j_k$ will become small, despite the fact that the unregularized problem may only imply that \textit{some} components of $\mathbf q$ cannot be stably recovered. To address this problem, we will appeal to a stochastic (Bayesian) interpretation of the inverse problem
\cite{Tar04,KS05}. From this perspective, we assume that our measurements $\mathbf z$ are stochastic because they are corrupted by noise, and that consequently our recovered coefficients $\mathbf q$ are also stochastic variables whose joint probability distribution we would like to infer. If we assume that the components of $\mathbf z$ are distributed according to $\mathbf z = \mathbf {\hat z} + N(0,\Sigma^2)$ with a nominal (but unknown) ``exact'' measurement value $\mathbf{\hat z}$ (that is, we assume that the noise is Gaussian and that our guessed noise levels $\sigma_\ell$ are indeed correct), then the desired probability distribution for $\mathbf q$ will be of the form
\begin{equation}
   \label{eq:finite-dimensional-def-prob}
    p(\mathbf q | \mathbf z) = \kappa e^{-{\mathcal{J}}_\text{red}(\mathbf q)}
\end{equation}
where $\kappa$ is a normalization constant whose concrete value is not
of importance to us.

Given this interpretation, the question of how much we know about the
individual components of $\mathbf q$ can be related to the
\textit{uncertainty} under $p(\mathbf q | \mathbf z)$ -- namely, we
should choose the information weights $j_k$ as the \textit{inverse of
  the standard deviation%
  \footnote{The inverse of the variance is often called the ``precision'' with which a parameter is known \cite{degroot1970}.}
  of $q_k$}, that is, equal to
$\frac{1}{\sqrt{\text{var}_p(\mathbf q)_k}}$, where
\begin{equation}
  \label{eq:def-jk-via-cov}
    \text{var}_p(\mathbf q)_k = 
    \int \left(q_k-{\mathbb E}_p[\mathbf q]_k\right)^2 p(\mathbf q| \mathbf z)  \;\text{d}q.
\end{equation}
Here, ${\mathbb E}_p[\mathbf q] = 
\int \mathbf q \; p(\mathbf q| \mathbf z) \;\text{d}q$.%
\footnote{Because we are considering a linear problem
\eqref{eq:finite-dimensional-state-eq} and because the objective
function ${\mathcal{J}}_\text{red}$ is quadratic, the expectation value
${\mathbb E}_p(q_k)$ in \eqref{eq:def-jk-via-cov} is equal to the
solution of the original, deterministic problem
\eqref{eq:finite-dimensional-solution}.}
This choice of $j_k$ has the pleasant property of making
\eqref{eq:finite-dimensional-bound} dimensionally correct also for
cases where the components of $\mathbf q$ have different physical
units. Yet, it turns out to be computationally difficult to obtain the variances 
$\text{var}_p(\mathbf q)_k$, and as a consequence we will show
  in the next section how we can define $j_k$ in a way that
  allows for an efficient computation.

\subsection{Estimating the information content for each component of the solution vector}
\label{sec:finite-dimensional-information}

The definition of an information content $j_k$ based on the inverse of
the variance in the stochastic inverse problem makes intuitive
sense. The question remains whether these weights $j_k$ can be
computed efficiently. As discussed in Remark~\ref{remark:why-not-cov} at the end of this
section, the answer is no, but we will show next that approximations can be found that are efficiently computable.

To do so, recall that the variances $\text{var}_p(\mathbf q)_k=\text{cov}_p(\mathbf q)_{kk}$ are the diagonal entries of the covariance matrix associated with $p(\mathbf q|\mathbf z)$, where the covariance matrix is defined as 
\begin{align}
    \label{eq:def-covariance-matrix}
    \text{cov}_p(\mathbf q)_{kl}
    =
    \int \left(q_k-{\mathbb E}_p[\mathbf q]_k\right)\left(q_l-{\mathbb E}_p[\mathbf q]_l\right) p(\mathbf q| \mathbf z) \;\text{d}q.
\end{align}

As mentioned above, in an ideal world, one would define an information density based on the variances $\text{var}_p(\mathbf q)_{k}$, for example by using the inverse of the standard deviation. However, the covariance matrix or a reasonably good approximation can generally not be computed in a computationally tractable way for large problems using sampling techniques such as MCMC unless one utilizes properties of the problem (such as when the forward model is linear and the objective function is quadratic). We do not want to base our approach on using properties that typical applications do not possess. Instead, we use techniques known from statistical analysis. In particular, one can estimate the covariance matrix via the Cram{\'e}r-Rao bound that states that 
\begin{gather*}
    \text{cov}_p(\mathbf q) \ge I_p^{-1}
\end{gather*}
in the sense that $[\text{cov}_p(\mathbf q) - I_p^{-1}]$ is a positive semidefinite matrix. Here, $I_p$ is the Fisher information matrix defined by
\begin{equation}
  \label{eq:fisher}
    (I_p)_{kl} = -{\mathbb E}\left[ \frac{\partial^2}{\partial q_k\partial q_l} \ln p(\mathbf q|\mathbf z) \right],
\end{equation}
which for our choice of $p(\mathbf q|\mathbf z)$ and ${\mathcal{J}}_\text{red}(\mathbf q)$ evaluates to
\begin{equation*}
    (I_p)_{kl}= {\mathbb E}\left[ \frac{\partial^2}{\partial q_k\partial q_l} \left( -\ln\kappa + {\mathcal{J}}_\text{red}(\mathbf q) \right) \right]
    = {\mathbb E}\left[ \frac{\partial^2}{\partial q_k\partial q_l} {\mathcal{J}}_\text{red}(\mathbf q) \right]
    = Q_{kl}.
\end{equation*}
In other words, the Fisher information matrix \textit{can} be computed efficiently, unlike the covariance matrix. At the same time, the estimate above requires us to compute the \textit{inverse} of the Fisher matrix, which for large and ill-conditioned problems is again not computable efficiently and accurately. However, we can use the the following inequality, again known from statistical analysis \cite{COV06}:
\begin{gather*}
    \text{var}_p(\mathbf q)_k = \text{cov}_p(\mathbf q)_{kk}
    \ge \left[I_p^{-1}\right]_{kk}
    \ge \left[(I_p)_{kk}\right]^{-1}.
\end{gather*}
These statements then provide us with an efficient way to estimate $\text{var}_p(\mathbf q)_k$:
\begin{equation}
\label{eq:ineq-variance-Q}
    \text{var}_p(\mathbf q)_k
    \ge \left[(I_p)_{kk}\right]^{-1}
    = \left[ Q_{kk} \right]^{-1}.
\end{equation}
In the spirit of the transition from \eqref{eq:information-bound} to
\eqref{eq:information-bound-approximation}, let us then define the
information content of the $k$th parameter as
\begin{equation}
\label{eq:finite-dimensional-information-content}
    j_k := \sqrt{Q_{kk}}.
\end{equation}

\begin{remark}
  \label{remark:compute-Q}
Based on the definition $Q=S^TA^{-T}M^T \Sigma^{-2} M A^{-1} S + \beta R^TR$,
the elements $j_k=\sqrt{Q_{kk}}$ can be computed in different ways by setting parentheses in the defining expression. The first way computes 
\begin{align*}
    Q_{kk} &= \mathbf e_k^T Q \mathbf e_k
    \\
    &=
    (\Sigma^{-1}M A^{-1} S \mathbf e_k)^T (\Sigma^{-1}M A^{-1} S \mathbf e_k) 
    + \beta (R \mathbf e_k)^T (R\mathbf e_k)
    \\
    &=
    (\Sigma^{-1}M A^{-1} \mathbf s_k)^T (\Sigma^{-1}M A^{-1} \mathbf s_k) 
    + \beta (R \mathbf e_k)^T (R\mathbf e_k)
    \\
    &=
    (\Sigma^{-1}M\mathbf h_k)^T (\Sigma^{-1}M\mathbf h_k) 
    + \beta \mathbf r_k^T \mathbf r_k
    \\
    &=
    \sum_\ell \frac{1}{\sigma_\ell^2} \left(\mathbf m_\ell^T\mathbf h_k\right)^2
    + \beta \mathbf r_k^T \mathbf r_k,
\end{align*}
where $\mathbf e_k$ is the $k$th unit vector and 
\begin{gather*}
 \mathbf h_k = A^{-1}S \mathbf e_k 
 = A^{-1} \mathbf s_k,
 \qquad\qquad
 \mathbf r_k = R\mathbf e_k.
\end{gather*}
That is, computing the information content vector $\mathbf j$ requires the solution of the forward operator $A$ for each of the source terms,
plus a few matrix vector products.

An alternative way involves computing $MA^{-1}=(A^{-T}M^T)^T$ first. Because the vectors $\mathbf m_\ell$ form the rows of $M$ (and so the columns of $M^T$), we can compute vectors
\begin{gather*}
 \mathbf h_\ell^\ast = A^{-T} \mathbf m_\ell,
\end{gather*}
and then recognize that
\begin{align*}
    Q_{kk}
    &=
    \sum_\ell \frac{1}{\sigma_\ell^2}\left((\mathbf h_\ell^\ast)^T \mathbf s_k\right)^2 
    + \beta \mathbf r_k^T \mathbf r_k.
\end{align*}
This approach requires solving a linear system with $A^T$ for each measurement.

Which of the two ways of computing $Q_{kk}$ is more efficient depends on whether there are more measurements than source terms, or the other way around.%
\footnote{Clearly, both ways are expensive for real-world cases with many parameters and many measurements. We will come back to this in our conclusions and outlook, Section~\ref{sec:conclusions}.}
\end{remark}

Regardless of the way $Q_{kk}$ (and consequently $j_k=\sqrt{Q_{kk}}$)
is computed, it can be interpreted as having contributions from all
measurements (through the sum over $\ell$) and from
regularization. The scalar product $\mathbf m_\ell^T\mathbf h_k$ can
be considered as the influence of the forward propagated sources
($\mathbf h_k$) on measurements. On the other hand, the equivalent
term $(\mathbf h_\ell^\ast)^T\mathbf s_k$ corresponds to a view where
we first compute an adjoint solution $\mathbf h_\ell^\ast$ that indicates which possible source terms affect a measurement functional, and then take the dot product with a concrete source $\mathbf s_k$. Both views represent the sensitivity of measurement functionals to sources.

Interestingly, the formula expresses the intuitive concept that \textit{information is additive}: If there are no measurements and no regularization, then $j=0$; each measurement in turn adds a non-negative contribution.
Finally, because the measurement contributions to $Q_{kk}$ are proportional to $\frac{1}{\sigma_\ell^2}$, and because the information content is the square root of $Q_{kk}$ we have the convenient and reasonable property that \textit{information is inversely proportional to the measurement uncertainty}. As a consequence, the definition of $j$ chosen here matches the goals we have set out at the end of the introduction for an information content or density.

\begin{remark}
\label{remark:why-not-cov}
    The \textit{intent} of the derivations above is to define the
    information content as $j_k \approx
    \frac{1}{\sqrt{\text{var}_p(\mathbf q)_k}}$, see the end of Section~\ref{sec:def-based-on-covariance}, but in practice, we define it as
    $j_k := \sqrt{Q_{kk}}$, the latter being an approximation of the former via the Cram{\'e}r-Rao bound and \eqref{eq:ineq-variance-Q}. 
    We do not choose the former because in practice computing all
    diagonal elements $Q_{kk}$ can be done efficiently via the
    techniques outlined in Remark~\ref{remark:compute-Q}, whereas
    computing all variances $\text{var}_p(\mathbf q)_k$ (i.e., all
    diagonal elements of the covariance matrix) is far more expensive,
    in particular for nonlinear problems (see also
    Appendix~\ref{sec:nonlinear}). Yet, it is conceivable that the variances can be efficiently \textit{approximated}, for example via randomized linear algebra techniques and sampling methods to compute expectation values. Whether it makes a \textit{practical difference} to use one or the other definition of $j_k$ -- and consequently whether it is worth the effort to actually compute the variances -- is, of course, a different question; we will leave the answer to this question to future research.
\end{remark}

\begin{remark}
  In practice, models are rarely linear. For nonlinear models, the Cram{\'e}r-Rao bound as well as \eqref{eq:ineq-variance-Q} still hold, but perhaps with a larger gap between left and right hand sides. Nonetheless, we can analyze the nonlinear case. Such an analysis is provided in Appendix~\ref{sec:nonlinear} in substantial detail, but let us preview the results here already: Using standard statistical and analytical techniques, the definition $j_k := \sqrt{Q_{kk}}$ we have chosen here will still be a good approximation of the inverse of the variance of parameter $q_k$ if (i) the model is not too nonlinear, or (ii) if the noise level is relatively small. Furthermore, the computation of information contents $j_k$ only requires \textit{linearized} solves of the forward or adjoint model, making the relative cost substantially cheaper than in the linear case because the solution of the inverse problem now requires the solution of \textit{nonlinear} forward problems.
\end{remark}

\section{Extension to infinite-dimensional inverse problems}
\label{sec:infinite-dimensional}

We can extend the reasoning of the previous section to infinite-dimensional inverse problems. Specifically, let us consider the linear source identification problem
\begin{align}
    \label{eq:infinite-dimensional-state-eq}
    {\mathcal{L}} u(\mathbf x) = \sum_k q_k s_k(\mathbf x), 
    \qquad\qquad \forall \mathbf x \in \Omega,
\end{align}
where $\mathcal{L}$ is a differential operator acting on functions defined on a domain $\Omega\subset \R^d$, and the equations are augmented by appropriate boundary conditions on $\partial\Omega$ whose details we will skip for the moment. As before, $s_k(\mathbf x)$ are possible source vectors and $q_k$ are their relative strengths. We again seek to identify source strengths $q_k$. Importantly, we will assume that the source terms are all of the form
\begin{align*}
    s_k (\mathbf x) = \chi_{\omega_k}(\mathbf x),
\end{align*}
where $\chi_{\omega_k}$ is the characteristic function of a subdomain $\omega_k$, and we assume that $\omega_k\cap\omega_l = \emptyset$ for $k\neq l$ and $\bigcup_k \overline\omega_k = \overline\Omega$. In other words, we seek to identify a source term that is a piecewise constant function defined on a partition of the domain $\Omega$.

We will infer the source strengths $q_k$ through
(linear, noisy) measurements
\begin{align}
    \label{eq:infinite-dimensional-measurements}
    z_\ell = 
    \left<m_\ell,u\right> + \varepsilon_\ell = 
    \int_\Omega m_\ell(\mathbf x) u(\mathbf x) \; \text{d}x + \varepsilon_\ell, 
\end{align}
where $\varepsilon_\ell$ is again measurement noise assumed to have
magnitude $\sigma_\ell$. The formalism we will develop will allow us
to assign an \textit{information content} to each $q_k$. Because the
source strengths $q_k$ correspond to characteristic functions $s_k$ of
subdomains $\omega_k$, the information content $j_k$ divided by the
volume $|\omega_k|$ will define an \textit{information density} $j(\mathbf x)$, for which we can consider the limit case $|\omega_k|\rightarrow 0$. This limit is not computable, but we can use finite subdivisions into regions $\omega_k$ that allow us to approximate it with reasonable accuracy.

To make these concepts concrete, in Section~\ref{sec:results} we will  consider this model where ${\mathcal{L}}$ is an advection-diffusion operator, ${\mathcal{L}} = -D\Delta + \mathbf b \cdot \nabla$, and where $z_\ell$ correspond to point measurements of $u(\mathbf x)$ (or a well-defined approximation of point measurements if the solution $u$ is not guaranteed to be a continuous function). This example is motivated by a desire to identify sources of air pollution from sparse measurements at a finite number of points.

\subsection{Definition of the inverse problem}
\label{sec:infinite-dimensional-definition}

The  inverse problem we have described above then has the following mathematical formulation, where we also include an $L_2$ regularization term:
\begin{align}
  \label{eq:inverse-problem}
\begin{split}
\underset{\mathbf q, u}{\text{min}}\; \mathcal J(\mathbf q, u) &= 
 \frac 12 \sum_\ell \frac{1}{\sigma_\ell^2} \left|\left< m_\ell, u \right> - z_\ell\right|^2
 + \frac{\beta}{2} \left\|\sum_k q_k s_k\right\|_{L_2(\Omega)}^2,
  \\
  \text{such that}\quad
  &
 {\mathcal{L}} u = \sum_k q_k s_k.
 \end{split}
\end{align}

 As in the previous section, we can define a reduced objective function
\begin{align}
    \label{eq:inverse-reduced-objective}
    {\mathcal{J}}_\text{red}(\mathbf q) &= 
    {\mathcal{J}}\left(\mathbf q, {\mathcal{L}}^{-1} \sum_k q_k s_k\right)
    \\
    \notag
    &=
    \frac 12 \sum_\ell \frac{1}{\sigma_\ell^2}\left|\left< m_\ell, {\mathcal{L}}^{-1} \sum_k q_k s_k \right> - z_\ell\right|^2
    + \frac{\beta}{2} \left\|\sum_k q_k s_k\right\|_{L_2(\Omega)}^2,
\end{align}
which gives rise to a related stochastic inverse problem with a probability density $p(\mathbf q|\mathbf z)$ defined as in \eqref{eq:finite-dimensional-def-prob}.

\subsection{Defining the information content}
\label{sec:derivation}

As in Section~\ref{sec:finite-dimensional-information}, we can again identify the \textit{information content} associated with each parameter $q_k$ via the precision, i.e., inverse of the variance $\text{var}_p(\mathbf q)_k = \text{cov}_p(\mathbf q)_{kk}$, and the estimate we have for the variance based on the Fisher information matrix. 

In the finite-dimensional case, the Fisher information matrix $I_p$ could be computed by solving one forward problem for each source vector $\mathbf s_k$. The same is true for the current infinite-dimensional situation:
\begin{proposition}
  For the model problem defined above, the Fisher information matrix $I_p$
  defined in \eqref{eq:fisher} has the following form:
  \begin{align}
    \label{eq:fisher-model}
    (I_p)_{kl}
    &=
    Q_{kl}
    \intertext{where}
    \label{eq:Q-matrix}
    Q_{kl} &=
    \sum_\ell \frac{1}{\sigma_\ell^2} \left< m_\ell, h_k \right>
            \left< m_\ell, h_l \right>
    + \beta \int_\Omega s_k s_l,
  \end{align}
  and where $h_k$ satisfies the equation
  \begin{align}
    \label{eq:def_h_k}
    {\mathcal{L}} h_k(\mathbf x) &= s_k(\mathbf x)
     \qquad\qquad \forall \mathbf x \in \Omega,
  \end{align}
  again augmented by appropriate boundary conditions for $h_k$.
\end{proposition}

{\parindent0pt
\textbf{Proof.}
Recall that
\begin{gather*}
    (I_p)_{kl} = -{\mathbb E}\left[ \frac{\partial^2}{\partial q_k\partial q_l} \ln p(\mathbf q|\mathbf z) \right],
    \qquad\qquad
    \text{with}\qquad
        p(\mathbf q | \mathbf z) = \kappa e^{-{\mathcal{J}}_\text{red}(\mathbf q)}.
\end{gather*}
Based on the definition of ${\mathcal{J}}_\text{red}(\mathbf q)$ and the linearity of $\mathcal{L}$, we then obtain
\begin{align*}
  (I_p)_{kl}
  &=
  \frac{\partial^2}{\partial q_k\partial q_l} {\mathcal{J}}(\mathbf q)
  \\
  &=
  \frac{\partial^2}{\partial q_k\partial q_l}
  \left[
  \frac 12 \sum_\ell \frac{1}{\sigma_\ell^2}\left< m_\ell, {\mathcal{L}}^{-1} \sum_r q_r s_r \right>^2
    + \frac{\beta}{2} \left\|\sum_r q_r s_r\right\|_{L_2(\Omega)}^2
    \right]
  \\
  &=
  \sum_\ell \frac{1}{\sigma_\ell^2}
            \left< m_\ell, {\mathcal{L}}^{-1} s_k \right>
            \left< m_\ell, {\mathcal{L}}^{-1} s_l \right>
    + \beta \int_\Omega s_k s_l,
\end{align*}
as claimed when using $h_k := {\mathcal{L}}^{-1} s_k$.
\hfill$\square$
}

\begin{proposition}
  The matrix $Q$ can alternatively be expressed through the following formula:
  \begin{align}
    \label{eq:Q-matrix-adjoint}
    Q_{kl} &=
    \sum_\ell \frac{1}{\sigma_\ell^2} \left< h_\ell^\ast, s_k \right>
            \left< h_\ell^\ast, s_l \right>
    + \beta \int_\Omega s_k s_l.
  \end{align}
  Here $h_\ell^\ast$ satisfies the equation
  \begin{align}
    \label{eq:def_h_k-adjoint}
    {\mathcal{L}}^\ast h_\ell^\ast(\mathbf x) &= m_\ell(\mathbf x)
     \qquad\qquad \forall \mathbf x \in \Omega,
  \end{align}
  where ${\mathcal{L}}^\ast$ is the adjoint operator to $\mathcal{L}$, using appropriate boundary conditions for $h_\ell^\ast$.
\end{proposition}

{\parindent0pt
\textbf{Proof.}
The proposition follows from the observation that
\begin{align*}
    \left< m_\ell, h_k \right>
    =
    \left< m_\ell, {\mathcal{L}}^{-1} s_k \right>
    =
    \left< {\mathcal{L}}^{-\ast} m_\ell, s_k \right>
    =
    \left< h^\ast_\ell, s_k \right>.
\end{align*}
\hfill$\square$
}

\begin{remark}
  As in the finite-dimensional case, the Fisher information matrix \eqref{eq:fisher-model} is easy to compute for problems with either not too many
  parameters (using \eqref{eq:Q-matrix}) or not too many measurements (then using \eqref{eq:Q-matrix-adjoint}). In either case, the functions $h_k$ or $h^\ast_\ell$ can be computed independently in parallel. Which of the two forms is more efficient depends on whether there are more measurements than source terms or the other way around. That said, in the discussions below, we will want to let $|\omega_k|\rightarrow 0$ and consequently make the number of source terms and parameters infinite, and in that case the adjoint formulation in \eqref{eq:Q-matrix-adjoint} provides the more useful perspective.
\end{remark}

The Fisher information matrix approximates the inverse of
the covariance matrix, and the diagonal elements of the Fisher matrix $I_p$
therefore provide a means to estimate the certainty in the corresponding
parameters $q_k$. In the same way as for the finite-dimensional case in 
\eqref{eq:finite-dimensional-information-content}, we can then define an information content for the parameter $k$ via
\begin{equation}
\label{eq:infinite-dimensional-information-content}
    j_k := \sqrt{Q_{kk}},
\end{equation}
where now
\begin{gather*}
    Q_{kk}
    = 
    \sum_\ell \frac{1}{\sigma_\ell^2}\left< m_\ell, h_k \right>^2
    + \beta \int_\Omega s_k^2
    =
    \sum_\ell \frac{1}{\sigma_\ell^2}\left< h^\ast_\ell, s_k \right>^2
    + \beta \int_\Omega s_k^2.
\end{gather*}

\subsection{Defining the information density}

The discussions in the previous section did not make use of any particular properties of the source basis functions $s_k$. We now examine the special case of indentifying a piecewise constant source function, i.e., where
\begin{align*}
    s_k (\mathbf x) = \chi_{\omega_k}(\mathbf x).
\end{align*}
In this case, the information content for the parameter $q_k$ associated with area $\omega_k$ is
\begin{equation}
\label{eq:infinite-dimensional-information-content-2}
    j_k
    =
    \sqrt{Q_{kk}}
    =
    \sqrt{\sum_\ell \frac{1}{\sigma_\ell^2}\left(\int_{\omega_k }h^\ast_\ell(\mathbf x) \; dx \right)^2
    + \beta |\omega_k|}.
\end{equation}
Since this quantity scales with the size of the subdomains $\omega_k$, it is reasonable to define a piecewise constant
\textit{information density} as
\begin{align}
\label{eq:infinite-dimensional-information-content-3}
  j(\mathbf x)|_{\omega_k} 
  &= \frac{1}{|\omega_k|} j_k 
  = \frac{1}{|\omega_k|} \sqrt{Q_{kk}}
  \notag
  \\
  &= \sqrt{\sum_\ell \frac{1}{\sigma_\ell^2}\left(\frac{1}{|\omega_k|} \int_{\omega_k }h^\ast_\ell(\mathbf x) \; dx \right)^2
    + \beta\frac{1}{|\omega_k|}}
  \\
  \label{eq:information-density-definition}
  &\approx
  \sqrt{\sum_\ell \frac{1}{\sigma_\ell^2} h_\ell^\ast(\mathbf x)^2 + \beta\frac{1}{|\omega_k|}}.
\end{align}

We can make a number of observations based on these definitions, analogous to the finite-dimensional case of the previous section:
\begin{itemize}
\item Because the definition of $h_\ell^\ast$ is independent of the choice of     $\omega_k$, the formulas shown above can be interpreted as saying that the information density has a component that results from the measurements $\ell$ (and, in particular, grows monotonically with the number of measurements), and a component that results from regularization.

\item As before, the information density is inversely proportional to the measurement uncertainties $\sigma_\ell$, in the absence of regularization.

\item Regularization bounds the amount of information from below:
  $j(\mathbf x)|_{\omega_k} \ge \sqrt{\beta/|\omega_k|}$. This dependence on
  the square root of the regularization parameter is well known \cite{Cox74}.

  \item The information content $j_k$ for a subdomain $\omega_k$ decreases under mesh refinement (i.e., as  the subdomains become smaller and smaller). This makes sense if we have only finitely many measurements, and the rather weak $L_2$ regularization term. In order to ensure the continued well-posedness (i.e., an information content $j_k$ that is bounded away from zero) under mesh refinement, we also need to either increase the number of measurements accordingly, or use a stronger regularization term; for the latter, we need to choose regularizers that are ``trace-class'' -- see also the discussion in the Conclusions (Section~\ref{sec:conclusions}).
\end{itemize}

\begin{remark}
\label{rem:information-is-independent-of-noise}
  An important observation is that the definitions of information
  content and information density above depend only on the forward
  operator and the measurement functionals, but not on 
  concrete measured values $z_\ell$. Consequently, and as anticipated, information
  quantities can be computed \textit{before} measurements, and they are independent of the specific noise in measurements later obtained. We will come back to this point in Sections~\ref{sec:results-def-inverse-problem} and \ref{sec:mesh-refinement-comparison},
  as well as in the discussion of nonlinear problems in Appendix~\ref{sec:nonlinear}.
\end{remark}

\begin{remark}
We remark that the idea of using variances
$Q_{kk}$ as a spatially variable measure of certainty is not new. For
example, \cite[Section 3.1]{Rawlinson_2014} illustrates the spatially
variable variance for a seismic imaging problem. Yet, the authors' definition is unclear regarding the role of regularization, misses the square root, and is then discarded as not very useful.
\end{remark}

\section{Using information densities}
\label{sec:examples-of-use}

Having shown a way to define an information density $j(\vec x)$, the
question is whether it is useful. Indeed, there are numerous
questions related to the practical solution of inverse problems for
which information densities could be useful. In the following
subsections, we therefore first outline three vignettes of
situations in which the information density could be useful. In Section~\ref{sec:results}, we then expand on one of these ideas using a concrete numerical example.

In the examples below, we will consider the situation in which we have discretized the source term $q(\mathbf x)=\sum_k q_k s_k(\mathbf x)$ in \eqref{eq:infinite-dimensional-state-eq} on a ``mesh'' $\mathbb T$, as is common in the finite element method. Because there are no differentiability requirements on $q(\mathbf x)$, it is common to identify the source term as a piecewise constant function on the mesh, and in this case, the source functions $s_k$ are the characteristic functions of the cells $K$ of $\mathbb T$. By identifying the index $k$ with a cell $K$, \eqref{eq:infinite-dimensional-information-content} and \eqref{eq:infinite-dimensional-information-content-2} then define an ``information content'' $j_K$ for each cell $K$ of the mesh.

\subsection{Using information densities for regularization}
\label{sec:vignette-regularization}
As a first example of where we believe that information densities could be used, we consider the regularization of inverse problems. In a large number of  practical applications, one regularizes inverse problems by adding a penalty term to the misfit function for the purpose of penalizing undesirable aspects of the recovered function. For example, in our definition of the source identification problem in Section~\ref{sec:infinite-dimensional-definition} (see also equation~\eqref{eq:inverse-problem}), we have penalized the \textit{magnitude} of the source term to be identified. The strength of this penalization is provided by the factor $\beta$.

A practical question is how large this factor $\beta$ should be. Many criteria have been proposed in the literature \cite{EHN96},
but, in practice many studies do not use any of these automatic criteria and instead choose values of $\beta$ that yield reasonable results based on trial and error.

Moreover, it is clear to many practitioners that regularization may
not be necessary to the same degree in all parts of the domain. For
example, if measurements are available only in parts of the domain
(say, on the boundary), then intuitively more information is
available to identify source strengths close to the boundary than deep
in the interior of the domain. A particularly obvious example is in
seismic imaging: There, we can accurately identify properties of
the Earth only in those places that are crossed by ray paths from
earthquake sources (predominantly located at plate boundaries) to
seismometer stations (predominantly located on land), but not in the
rest of the Earth
\cite{Auer2014,SambridgeRawlinson2005,Liu_2012,Rawlinson_2014,Bonadio_2021}. In
the definition of the information density $j(\mathbf x)$ in
\eqref{eq:information-density-definition}, this would imply that the
first term under the square root would be large along these ray paths,
but small elsewhere.  In cases such as this, a reasonable approach
would be to make the regularization parameter spatially variable:
large where little information is available, and small where
regularization is not as important, for example so that
$j(\mathbf x)\ge j_0$. This spatial variability could be achieved by replacing the
regularization term in \eqref{eq:inverse-problem} by
\begin{align*}
    \frac{1}{2} \left\|\sqrt{\beta(\mathbf x)} \left(\sum_k q_k s_k(\mathbf x)\right)\right\|_{L_2(\Omega)}^2,
\end{align*}
and defining $\beta(\mathbf x)$ in some appropriate way.

Using spatially variable regularization is not a new idea (see, for example, \cite{Auer2014,PMPOP99,ArridgeSchweiger93,THA10,WBBS05,YHL10,CGG97,Bonadio_2021}), although we are not aware of any references that would provide an overarching, systematic framework for choosing $\beta(\mathbf x)$.
In contrast, the connection between  information density $j(\mathbf x)$ and $\beta(\mathbf x)$ in \eqref{eq:information-density-definition} has the potential to provide such a systematic approach. A scheme based on this observation also satisfies other considerations that appear reasonable. For example, increasing the number of measurements, or decreasing the measurement error, leads to a larger information density and therefore to a smaller regularization term to satisfy $j(\mathbf x)\ge j_0$.

\subsection{Using information densities to guide the discretization of an inverse problem}
\label{sec:vignette-mesh-refinement}

In actual practice, inverse problems are solved by discretization. In our derivation above, we have chosen finitely many source functions $s_k$ that we have assumed are the characteristic functions of ``cells'' $K$ of some kind of mesh or subdivision of the domain $\Omega$ on which $q(\mathbf x)$ is defined, and then expanded the function we seek as
\begin{align*}
    q(\mathbf x) = \sum_k q_k s_k(\mathbf x).
\end{align*}

A practical question is how to create this subdivision. Oftentimes, the subdivision is chosen fine enough to resolve the features of interest but coarse enough to keep the computational cost in check. Regularization is frequently used to ensure that an overly fine mesh does not lead to unwanted oscillations in the recovered coefficients -- in other words, to keep the problem reasonably well-posed. 

Most often in the literature, the mesh for the inverse problem is either uniform or at least chosen a priori through insight into the
problem (for approaches in the latter direction, see for example, \cite{BDM10,Auer2014,Liu_2012}).
On the other hand, discretization is a form of regularization, and
it is reasonable to choose the mesh finer where more information is available -- say,
close to a measurement device -- but coarser where our measurements
have little information to offer. This idea has been used as a
heuristic in the past \cite{Auer2014}, or at least mentioned (see
Section 2.3 of \cite{Rawlinson_2014} and the references therein), but,
as with regularization, no overarching scheme is available to guide this choice of the mesh. \cite{Bigoni2020} is also an example where the mesh is made part of what needs to be estimated in a Bayesian inversion scheme, which in practice appears to lead to meshes that are more refined where information is available.
However, \cite{Bigoni2020} is concerned with choosing the mesh used for the solution of the \textit{state equation}, not the discretization of the parameter we seek -- although it seems reasonable to assume that the scheme could be adapted to the latter as well. At the same time, the scheme described in \cite{Bigoni2020} requires the solution of a Bayesian inverse problem, which is far more expensive than the deterministic approach we use here.

Moreover, an information density can provide such a guide to
determine optimal cell sizes. First, we conjecture that meshes
should be graded in such a way that the information content of each
cell (i.e., roughly the information density times the measure of a
cell) is approximately equal among all cells.  We will explore in
detail how well this works in practice in Section~\ref{sec:results}.
Second, in mathematical research, mesh refinement cycles are
frequently terminated whenever we run out of memory, out of patience,
or both, whereas in applications, mesh refinement is stopped whenever an
expert deems the solution sufficiently accurate. Either approach is
unsatisfactory, and the amount of information available per cell might
provide a more rational criterion to stop mesh refinement.

\subsection{Using information densities for experimental design}
\label{sec:vignette-experimental-design}

As a final example, we consider optimal experimental design, that is, the question of what, how, or where to measure so as to minimize the uncertainty in recovered parameters given a certain noise level in measurements.

Optimal experimental design for inverse problems is more difficult than for finite-dimensional parameter estimation problems because it is not entirely clear what the objective function should be when minimizing or maximizing by varying the specifics of measuring. For finite-dimensional problems, objective functions include the $A$-, $C$-, $D$-, $E$-, and $T$-optimality criteria, plus many variations \cite{AD92}. 

For the infinite-dimensional case (or discretized versions thereof),
the choice might be to maximize the information content in all of $\Omega$, or a subset $\omega\subset\Omega$:
\begin{align*}
    \phi_\text{information}(\{m_\ell\}) =
    \int_\omega j(\mathbf x) \; dx,
\end{align*}
where $\{m_\ell\}$ denotes the set of measurements to be performed and optimization will typically happen over a set of implementable such measurements.

If $\omega=\Omega$, then the integral above reduces to
\begin{align*}
    \phi_\text{information}(\{m_\ell\})
    = \sum_k \sqrt{Q_{kk}},
\end{align*}
based on the definitions in \eqref{eq:infinite-dimensional-information-content-2} and \eqref{eq:infinite-dimensional-information-content-3}. Recalling that the matrix $Q$ is the Fisher information matrix, we recognize that the criterion $\phi$ above is similar to -- but distinct from -- the generally not very frequently used $T$-optimality criterion that maximizes the sum of diagonal entries of $Q$ (i.e., the \textit{trace} of $Q$):
\begin{align*}
    \phi_T(\{m_\ell\})
    = \sum_k Q_{kk}.
\end{align*}

\section{Numerical examples of using information densities}
\label{sec:results}

The previous section provided three vignettes of how we imagine
information densities could be used for practical
computations. Exploring all of these ideas through numerical examples exceeds the reasonable length of a single publication, and as a consequence we will focus on only mesh refinement.

In the following subsections, we will first lay out the inverse problem we will use as a test case. We will then show numerical results that illustrate the use of information densities as applied to this problem.

All numerical results were obtained with a program that is based on
the open-source finite element library \texttt{deal.II} \cite{arndt2019dealii,Arndt2021}. This program is  available under an open source license as part of the \texttt{deal.II} code gallery at \url{https://dealii.org/developer/doxygen/deal.II/CodeGallery.html} under the name ``Information density-based mesh refinement''.

\subsection{The test case}

Let us consider the following question: Given an advection-diffusion problem for a concentration $u(\mathbf x)$, can we identify the sources $q(\mathbf x)$ of the concentration field from point measurements of $u$ at points $\mathbf \xi_\ell$? This kind of problem is widely considered in environmental monitoring of pollution sources \cite{Skaggs1994,Neupauer2000,Vrugt2008,Jamshidi2020}, and also when trying to identify the sources of nuclear radiation.  

Mathematically, we assume that the concentration field satisfies the stationary advection-diffusion equation
\begin{align}
\label{eq:forward-problem}
{\mathcal{L}} u(\mathbf x) \equiv
    \mathbf b(\mathbf x) \cdot \nabla u(\mathbf x)
    - D \Delta u(\mathbf x) &= q(\mathbf x)
    \qquad\qquad&& \text{in $\Omega$},
    \\
\intertext{where $\mathbf b$ is a (known) wind field and $D$ is the (known) diffusion constant. For simplicity, we will assume homogenous Dirichlet boundary conditions}
  u &= 0 && \text{on $\partial\Omega$}.
\end{align}

Concretely, for our computations, we will assume that $\Omega=(-1,1)^2\subset{\mathbb R}^2$ is a square, $D=1$, and $\mathbf b=(100,0)^T$. These choices lead to a P{\'e}clet number of 200; that is, the problem is advection dominated.

For the inverse problem, we ask whether we can recover the function
$q(\mathbf x)$ (or a discretized version of it) from measurements
$z_\ell$ at a number of points $\mathbf \xi_\ell \in \Omega,
\ell=1,\ldots,L$. That is, we consider
\eqref{eq:infinite-dimensional-measurements} with $m_\ell(\mathbf x) =
\delta(\mathbf x-\mathbf \xi_\ell)$ where
$\sigma_\ell$ is the assumed noise level for the measurement at
location $\xi_\ell$ (see below). We choose these points $\mathbf \xi_\ell$ equally distributed around two concentric circles of radius 0.2 and 0.6, centered at the origin, with 50 points on each of the circles, for a total of $L=100$ measurement points.

For our experiments, we will consider a situation where the data we
have, $z_\ell$, has been obtained by solving the forward problem with
the finite element method, using a known source distribution
$q^\ast(\mathbf x)$ that is equal to one in a circle of radius 0.2
centered at $(-0.25,0)^T$.  A solution $u^\ast(\mathbf x)$ can then be
evaluated at the points $\mathbf \xi_\ell$ to obtain ``synthetic'' measurements $z_\ell$ via
\begin{align}
\label{eq:synthetic-measurements}
    z_\ell = \left<m_\ell,u^\ast\right> + \varepsilon_\ell = u^\ast(\mathbf \xi_\ell) + \varepsilon_\ell.
\end{align}
We choose Gaussian noise $\varepsilon_\ell = N(0,\sigma_\ell^2)$ and set $\sigma_\ell=0.1 \max_{\mathbf x\in \Omega}|u^\ast(\mathbf x)|$.

\begin{figure}
    \centering
    \phantom.
    \hfill
    \includegraphics[height=0.33\textwidth]{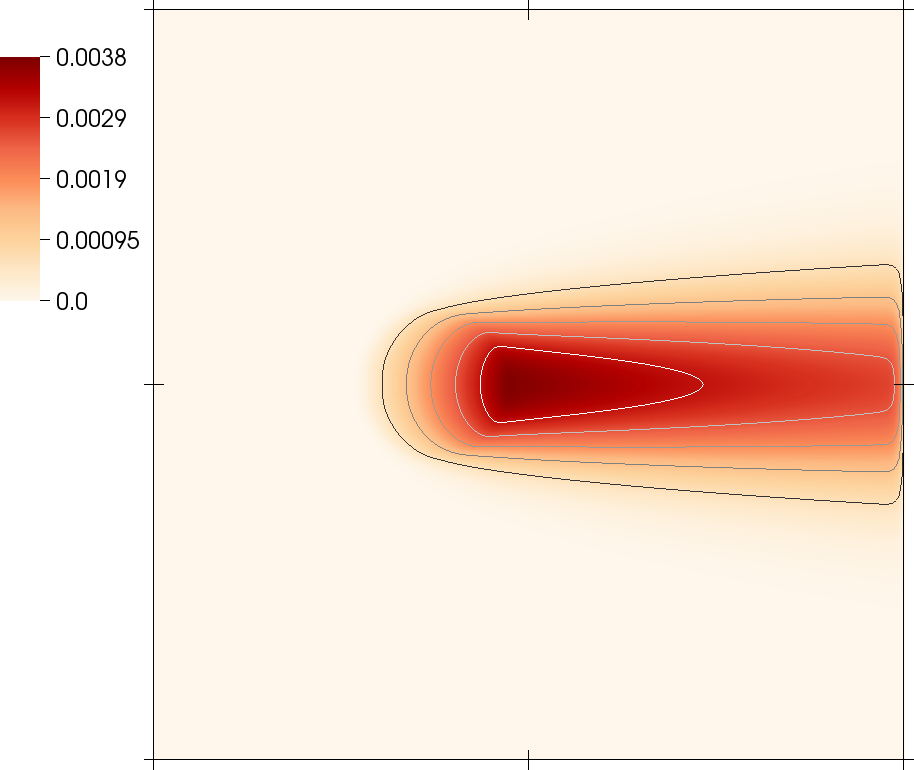}
    \hfill
    \includegraphics[height=0.33\textwidth]{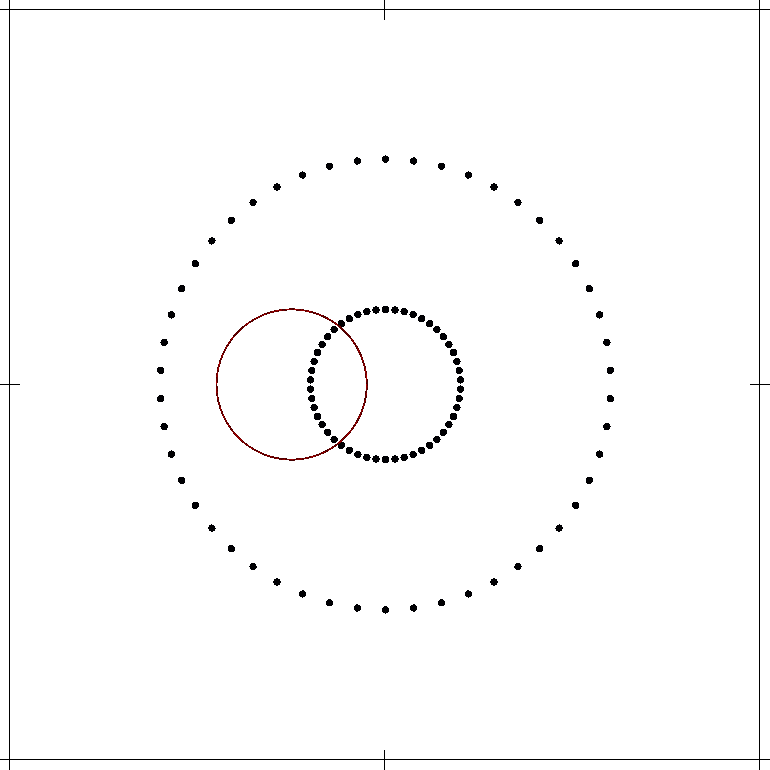}
    \hfill
    \phantom.
    \caption{\it Left: The solution $u^\ast(\mathbf x)$ of the forward problem from which we generate ``synthetic'' measurements $z_\ell$ via \eqref{eq:synthetic-measurements}. Right: The source
    term $q^\ast(\mathbf x)$ from which we compute synthetic measurements
    is constant and nonzero only in the solid red circle offset from the
    center; the detector locations $\xi_\ell, \ell=1,\ldots,L=100$ are
    marked by dots.}
    \label{fig:synthetic}
\end{figure}

To avoid an inverse crime, we solve for $u^\ast$ on a mesh that is different from the meshes used for all other computations. The solution of this forward problem so computed to obtain synthetic measurements is shown in Fig.~\ref{fig:synthetic}, along with the locations of the source term and the detector locations.

\subsection{The inverse problem}
\label{sec:results-def-inverse-problem}

The inverse problem we seek to solve is the identification of the source term $q(\mathbf x)$ (which we approximate via a finite-dimensional expansion $\sum_k q_k s_k(\mathbf x)$) in \eqref{eq:forward-problem}, based on the measurements described by \eqref{eq:synthetic-measurements}. We approach this problem by reformulating it in the form of the constrained optimization problem \eqref{eq:inverse-problem}, where we set the regularization parameter to $\beta=10^4$. This problem is then solved by introducing a Lagrangian
\begin{align*}
    {\mathfrak L}(u,q,\lambda)
    = {\mathcal{J}}(u,q) + \int_\Omega \lambda(\mathbf x)\left({\mathcal{L}}u(\mathbf x) - q(\mathbf x)\right) \text{d}x,
\end{align*}
and then solving the linear system of partial differential equations that results by setting the derivatives of ${\mathfrak L}$ to zero (that is, the optimality conditions). In strong form, these optimality conditions read
\begin{align}
\label{eq:inverse-problem-system}
\begin{split}
    {\mathcal{L}} u(\mathbf x) &= q(\mathbf x), \\
    {\mathcal{L}}^\ast \lambda(\mathbf x) &= -\sum_\ell \frac{1}{\sigma^2_\ell} (u(\mathbf \xi_\ell)-z_\ell) \delta(\mathbf x-\mathbf \xi_\ell), \\
    \beta q - \lambda &= 0.
    \end{split}
\end{align}
The solution of this system of equations is facilitated by discretizing on a finite element mesh. We use continuous, piecewise cubic elements for $u$ and $\lambda$, and discontinuous, piecewise constants elements for $q$.%
\footnote{This choice of higher order finite element spaces for $u$ and $\lambda$ is akin to solving for the forward and adjoint variables on a finer mesh than the source terms we seek to identify. As a result, we need not worry about satisfying discrete stability properties for the resulting saddle point problem. We can also, in essence, consider the forward and adjoint equation to be solved nearly exactly, with the majority of the discretization error resulting from the discretization of the source term $q(\mathbf x)$.}

\begin{figure}
    \centering
    \phantom{.}
    \hfill
    \raisebox{2cm}{$u(\mathbf x)$}
    \includegraphics[height=0.27\textwidth]{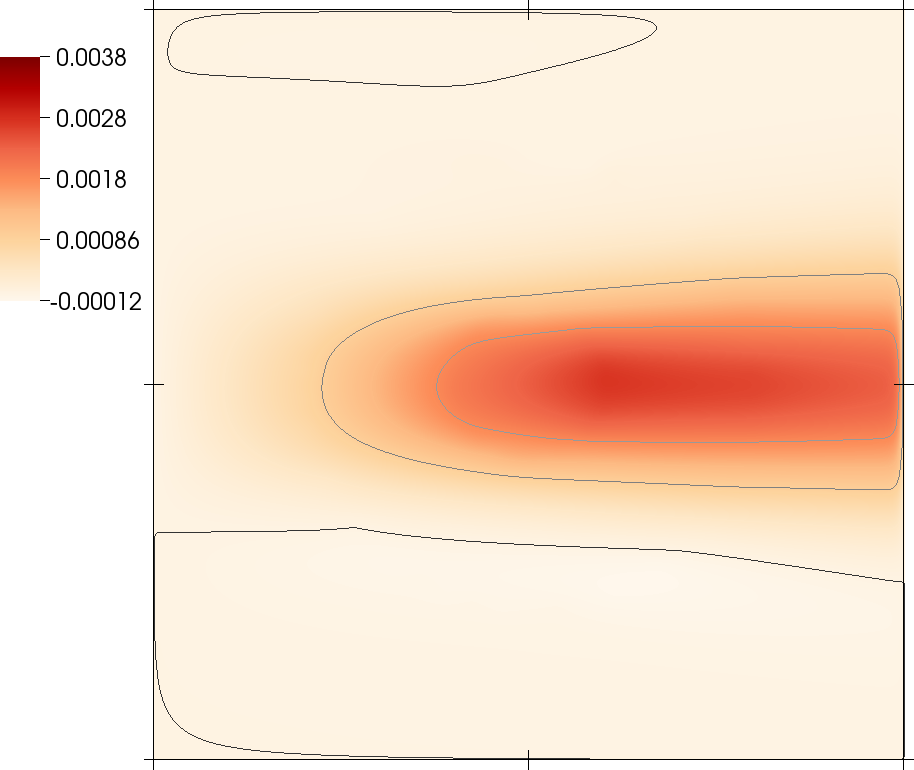}
    \hfill
    \includegraphics[height=0.27\textwidth]{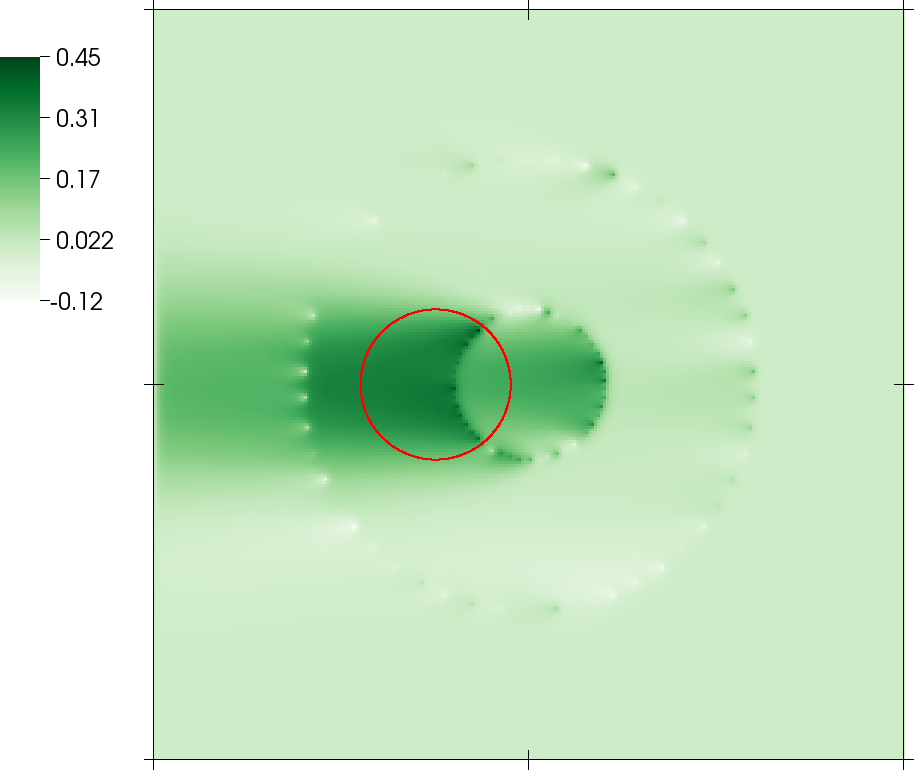}
    \raisebox{2cm}{$q(\mathbf x)$}
    \hfill
    \phantom{.}
    \\
    \phantom{.}
    \hfill
    \raisebox{2cm}{$\lambda(\mathbf x)$}
    \includegraphics[height=0.27\textwidth]{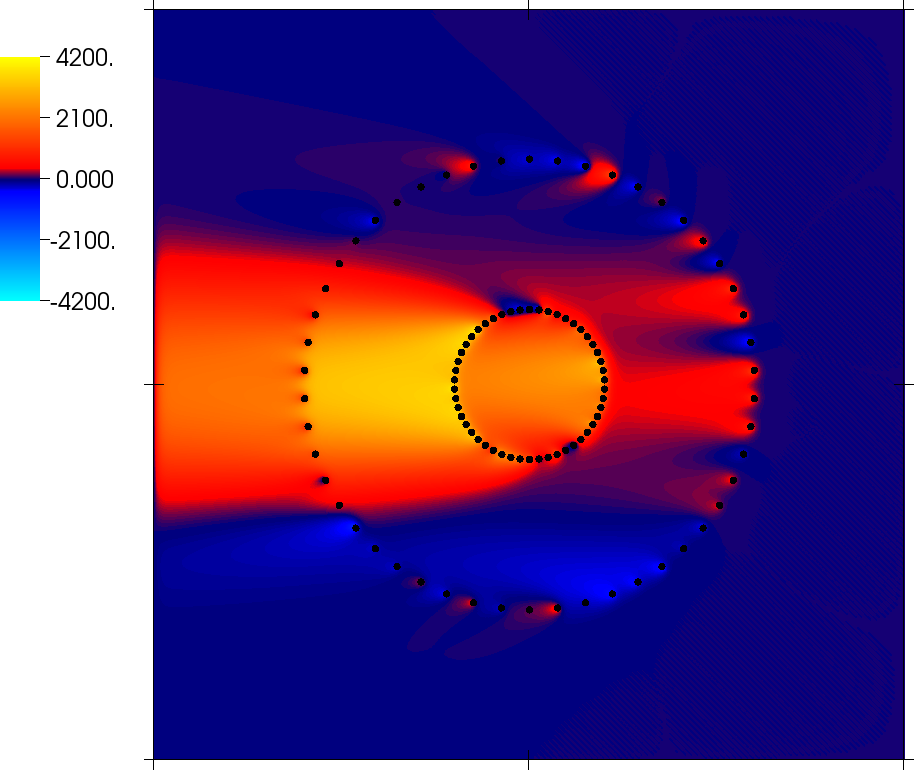}
    \hfill
    \includegraphics[height=0.27\textwidth]{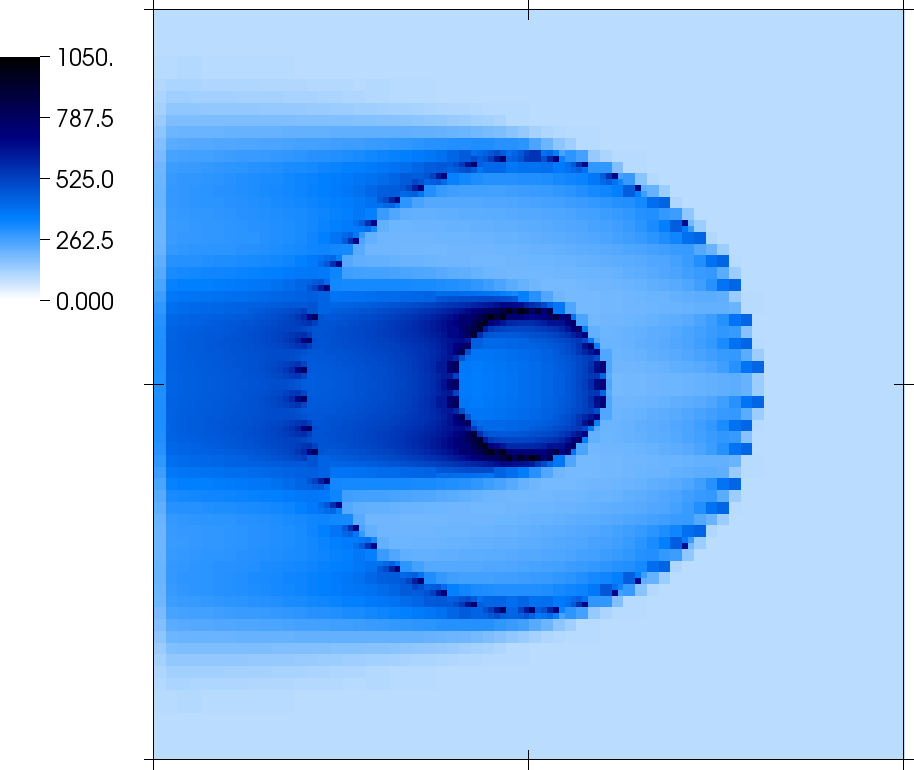}
    \raisebox{2cm}{$j(\mathbf x)$}
    \hfill
    \phantom{.}
    \caption{\it The solution of problem
    \eqref{eq:inverse-problem-system}, computed on a very fine finite element mesh. Top left: The primal variable $u(\mathbf x)$, shown with the same scale for color and isocontours as in Fig.~\ref{fig:synthetic}. Top right: The recovered sources $q(\mathbf x)$, i.e., the solution of the inverse problem. The red circle indicates the location of the source term used in generating the synthetic data. Bottom left: The adjoint variable $\lambda(\mathbf x)$. Bottom right: The information density $j(\mathbf x)$ associated with this problem, as defined in \eqref{eq:information-density-definition}.}
    \label{fig:inverse-solution-fine}
\end{figure}

The three components $u,q,\lambda$ of this solution, computed on a
very fine mesh with $256\times 256=\num{65536}$ cells, for which the
coupled problem has \num{1248258} unknowns, are shown in
Fig.~\ref{fig:inverse-solution-fine}. The maximal value of the
recovered source is less than half the maximal size of the ``true''
source, owing to the effect of the $L_2$ regularization term. As a
consequence, the forward solution $u$ is also too small. Furthermore,
the inverse problem places the source in a broader region than where
it really is, but this is not surprising: In an advection-dominated
problem, it is only possible to say with accuracy that the source is
\textit{upstream} of a detector, but not where in the upstream region
it actually is unless another detector further upstream indicates that
it must be downstream from the latter.

The adjoint variable $\lambda$ clearly illustrates the effect that the adjoint operator ${\mathcal{L}}^\ast$ transports information in the opposite direction $-\mathbf b$ of the forward operator ${\mathcal{L}}$, and that the sources of the adjoint equation are the residuals $-(u(\mathbf \xi_\ell)-z_\ell) \delta(\mathbf x-\mathbf \xi_\ell)$; here $u(\mathbf \xi_\ell)-z_\ell$ reflects the measurement error and, based on our choice of noise above, is Gaussian distributed with both positive and negative values.

The bottom right panel of the figure also shows the information density $j(\mathbf x)$ that corresponds to this problem, as defined in \eqref{eq:information-density-definition}. It illustrates that, given the location of detectors and the nature of the equation, information is primarily available upstream of detector locations. Notably, and as mentioned in Remark~\ref{rem:information-is-independent-of-noise}, the information density is based solely on the operator $\mathcal{L}$ and the measurement functionals $m_\ell$, but not on the actual measurements $z_\ell$ (or the noise that is part of $z_\ell$).

\subsection{Choice of mesh for the inverse problem}
\label{sec:results-mesh}

\begin{figure}
    \centering
    \includegraphics[width=0.19\textwidth]{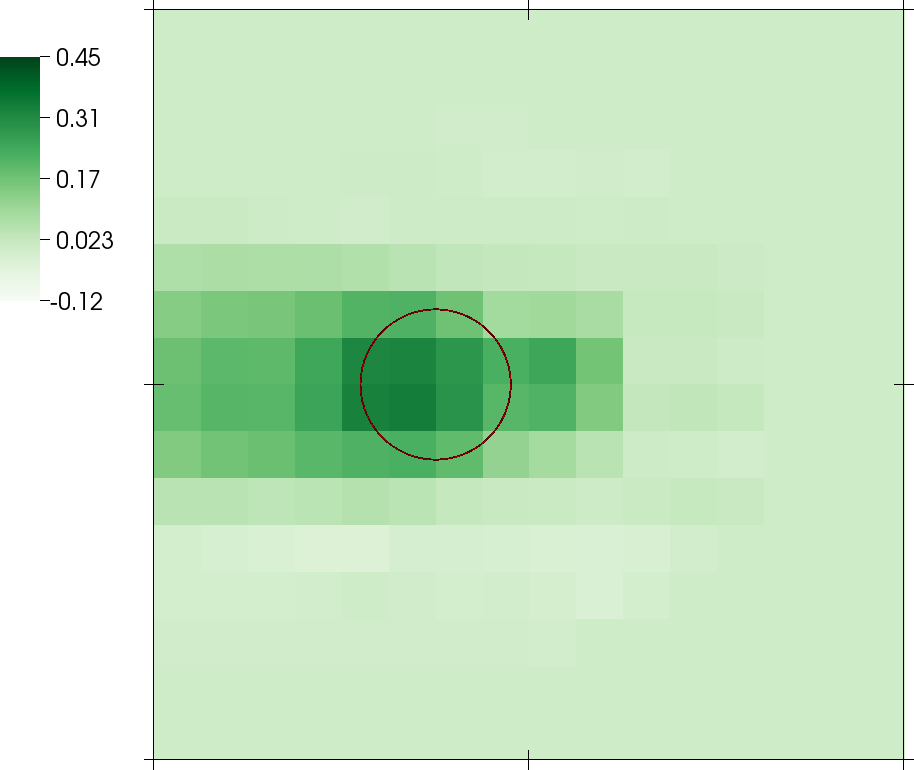}
    \hfill
    \includegraphics[width=0.19\textwidth]{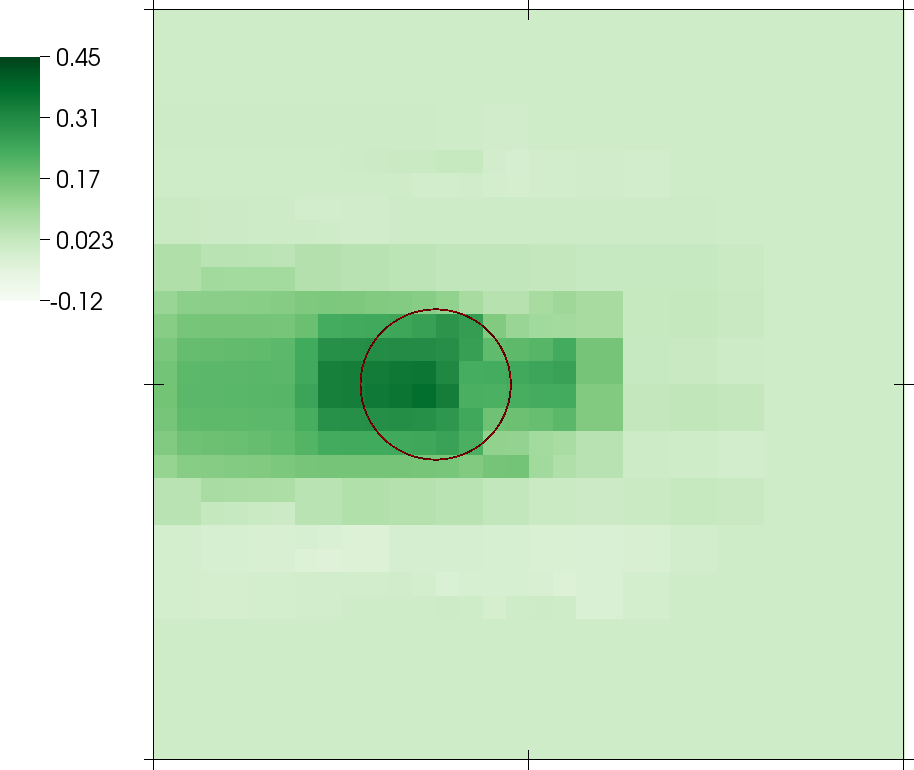}
    \hfill
    \includegraphics[width=0.19\textwidth]{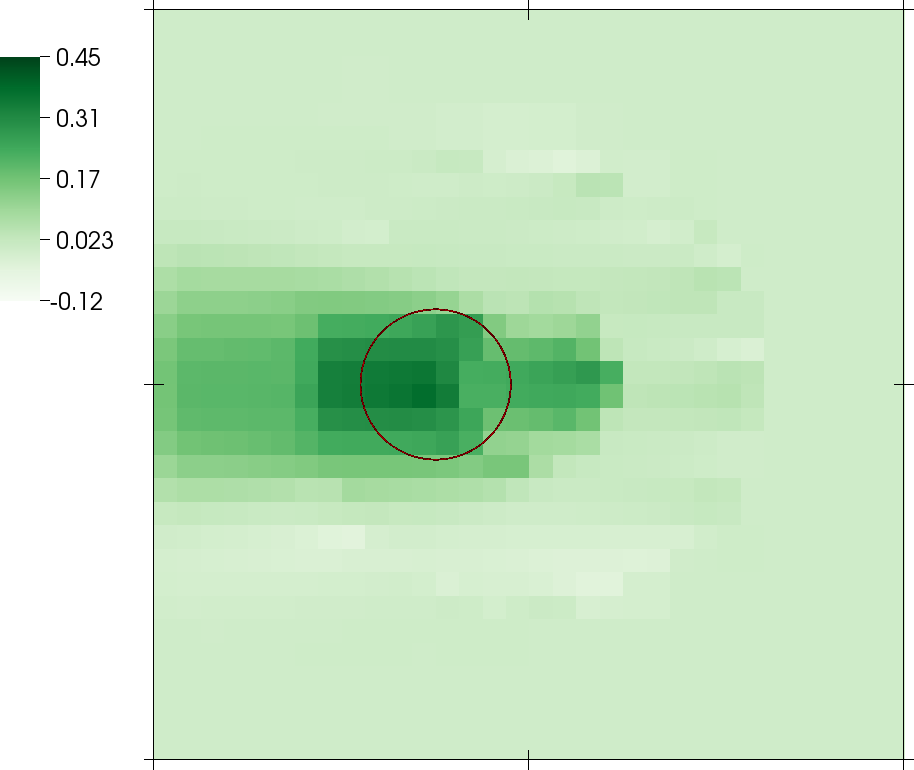}
    \hfill
    \includegraphics[width=0.19\textwidth]{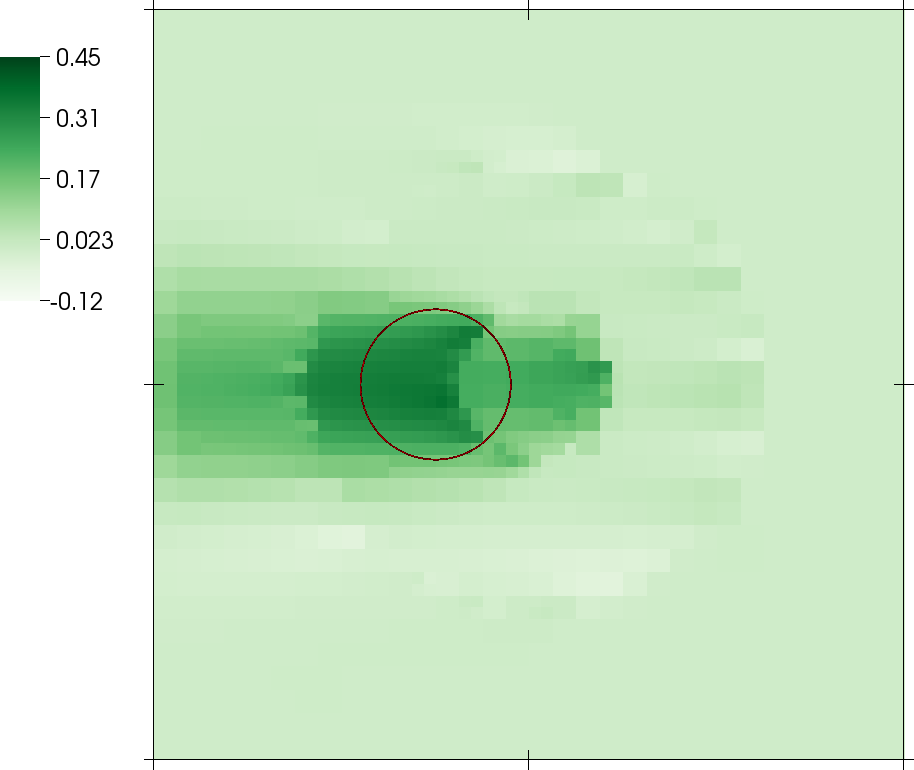}
    \hfill
    \includegraphics[width=0.19\textwidth]{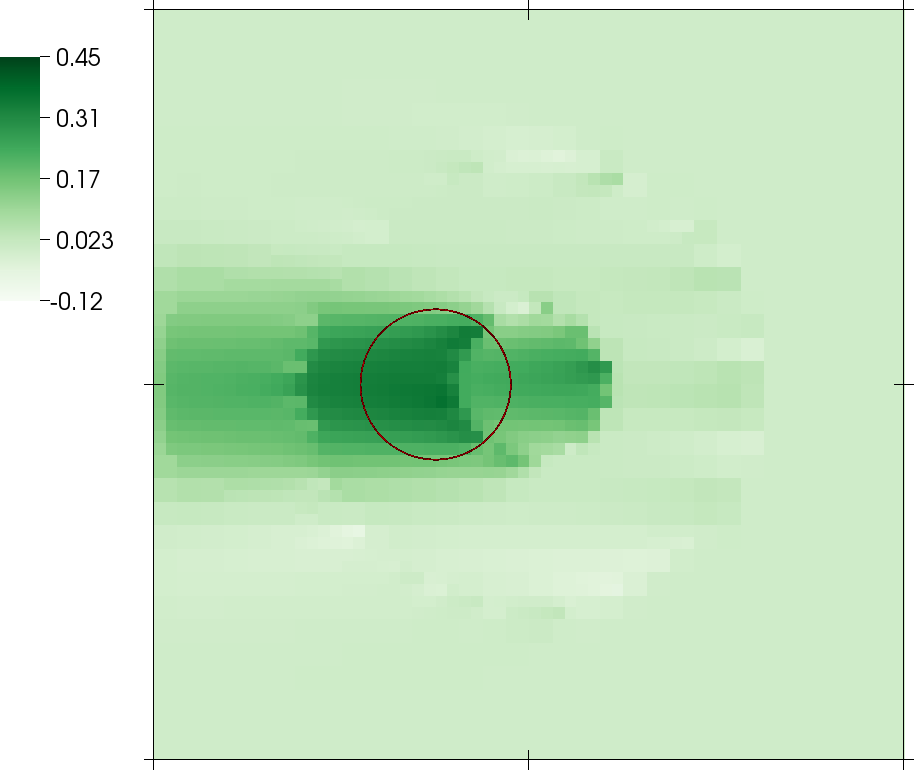}
    \\[12pt]
    \phantom.
    \hspace*{0.005\textwidth}
    \includegraphics[width=0.16\textwidth]{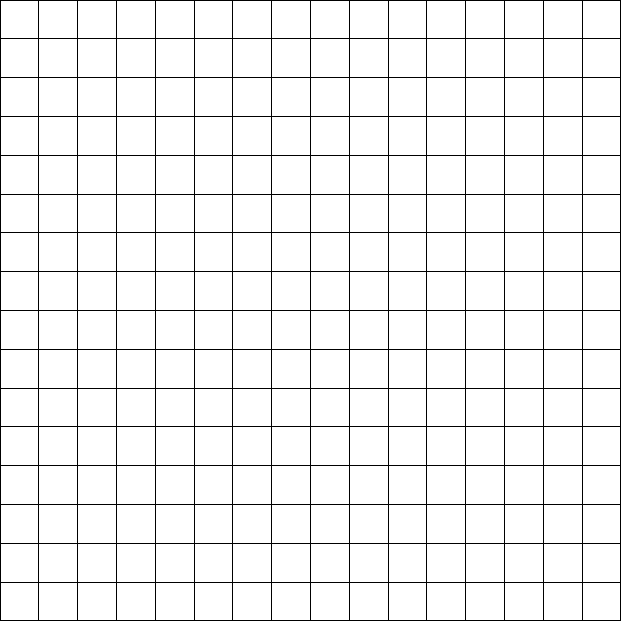}
    \hfill
    \includegraphics[width=0.16\textwidth]{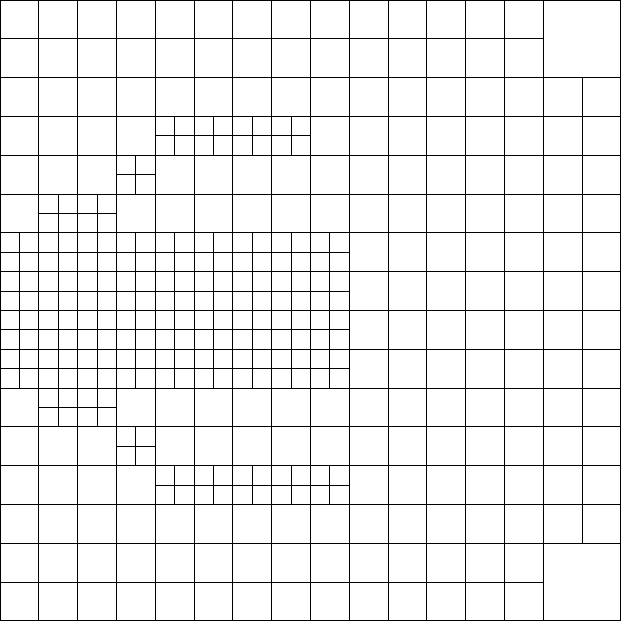}
    \hfill
    \includegraphics[width=0.16\textwidth]{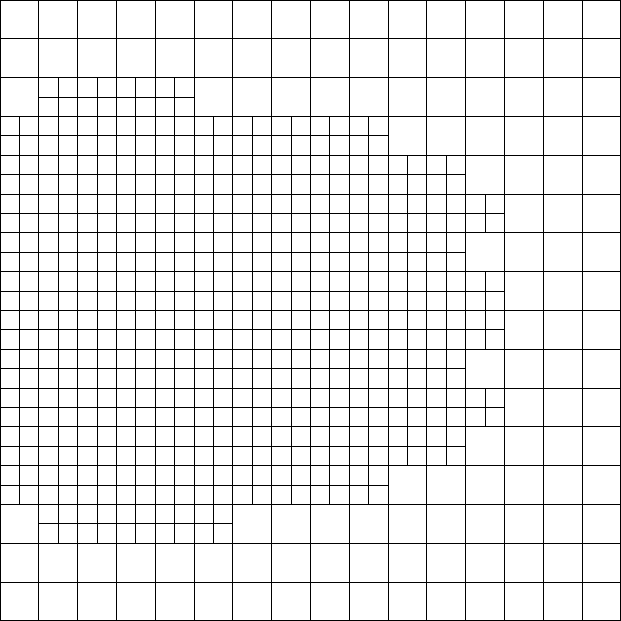}
    \hfill
    \includegraphics[width=0.16\textwidth]{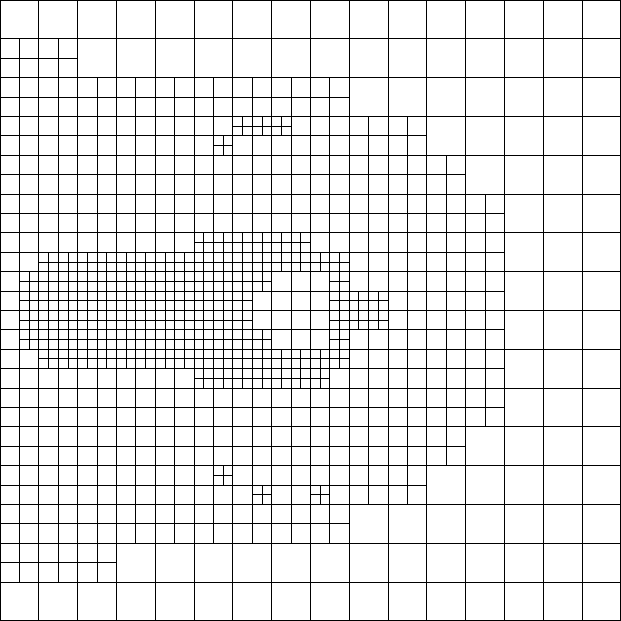}
    \hfill
    \includegraphics[width=0.16\textwidth]{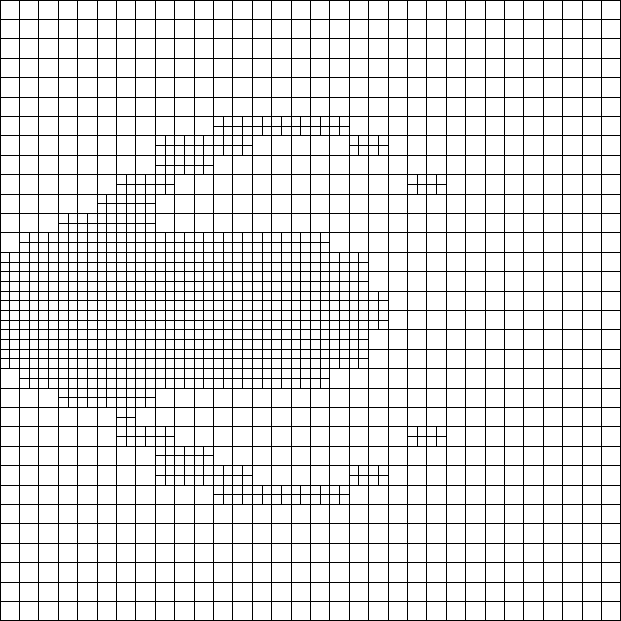}
    \caption{\it Reconstructions (top) on a sequence of meshes (bottom) refined based on the information content of each cell of the mesh.}
    \label{fig:inverse-solution-information-content-mesh}
\end{figure}

The question of interest then is how we can use information densities for mesh refinement. To answer this question, we have repeated the computations discussed above, but instead of using a uniformly refined mesh, we have used a sequence of meshes in which we refine cells hierarchically so as to equilibrate the information content $j_k$ of each cell $\omega_k$, see \eqref{eq:infinite-dimensional-information-content}, by always refining those cells that have the largest information content. The reconstructions and the sequence of meshes they are computed on are shown in Fig.~\ref{fig:inverse-solution-information-content-mesh}. For comparison with the computations mentioned above and shown in Fig.~\ref{fig:inverse-solution-fine}, the rightmost mesh has \num{1642} cells and the coupled problem solved on it has \num{32582} unknowns.

\subsubsection{Comparison with other mesh refinement criteria}
\label{sec:mesh-refinement-comparison}

The relevant question to ask is whether this mesh is better suited to the task than any other potential mesh. Answering this question is notoriously difficult in inverse problems because, in general, the exact solution of the problem is unknown if only finitely many measurements are available and if regularization is used. As a consequence, it is difficult to answer the question through comparison of convergence rates of different methods, for example.

\begin{figure}
    \centering
    \phantom.
    \hspace*{0.005\textwidth}
    \includegraphics[width=0.16\textwidth]{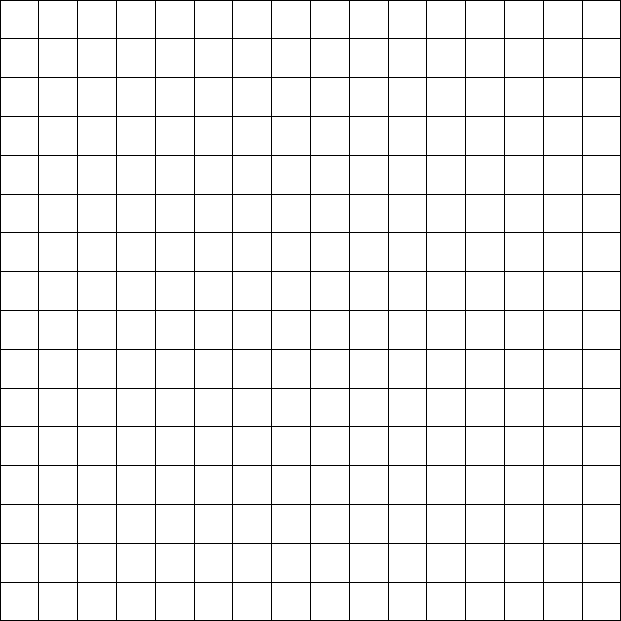}
    \hfill
    \includegraphics[width=0.16\textwidth]{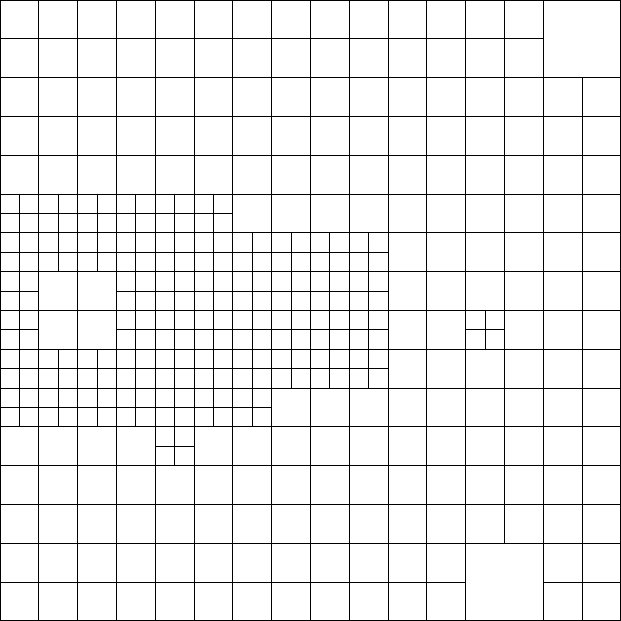}
    \hfill
    \includegraphics[width=0.16\textwidth]{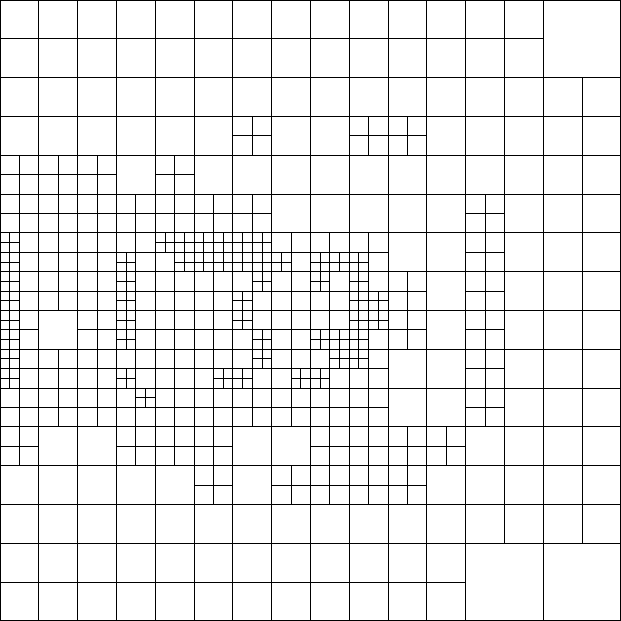}
    \hfill
    \includegraphics[width=0.16\textwidth]{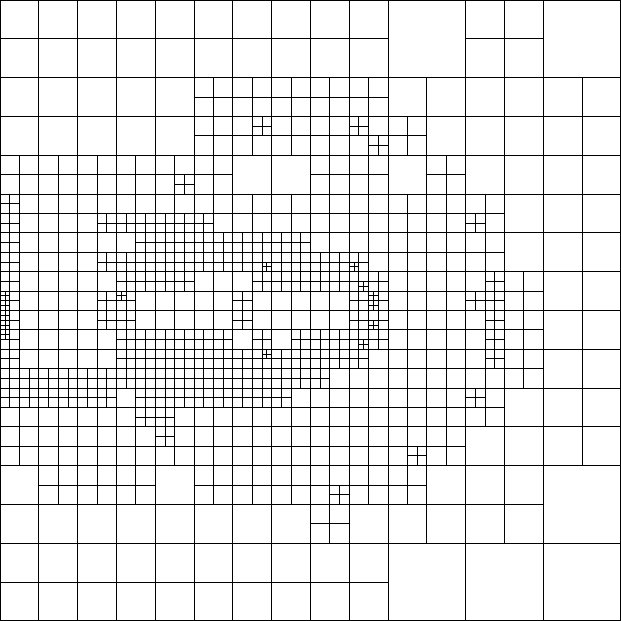}
    \hfill
    \includegraphics[width=0.16\textwidth]{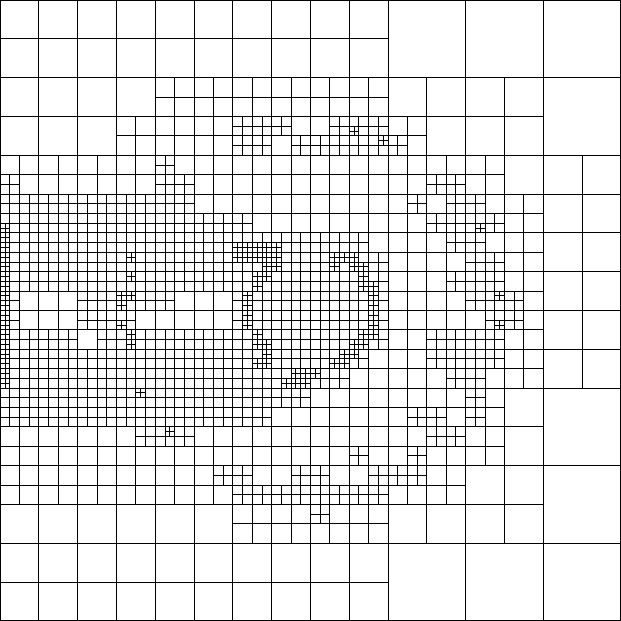}
    \\[12pt]
    \phantom.
    \hspace*{0.005\textwidth}
    \includegraphics[width=0.16\textwidth]{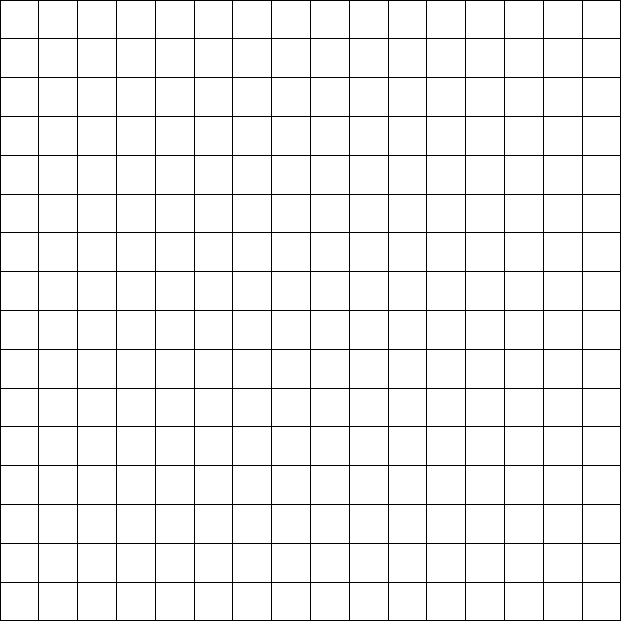}
    \hfill
    \includegraphics[width=0.16\textwidth]{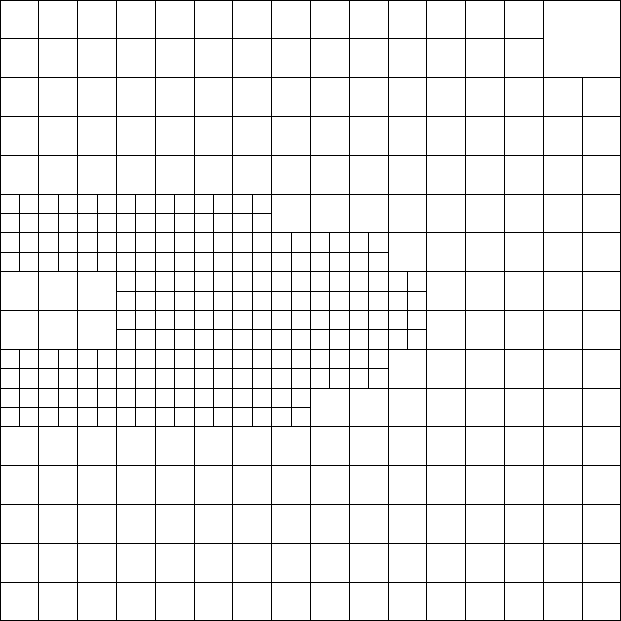}
    \hfill
    \includegraphics[width=0.16\textwidth]{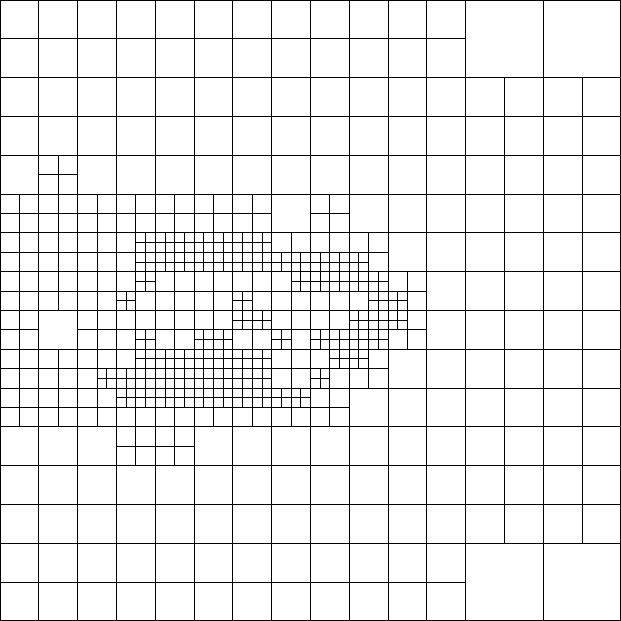}
    \hfill
    \includegraphics[width=0.16\textwidth]{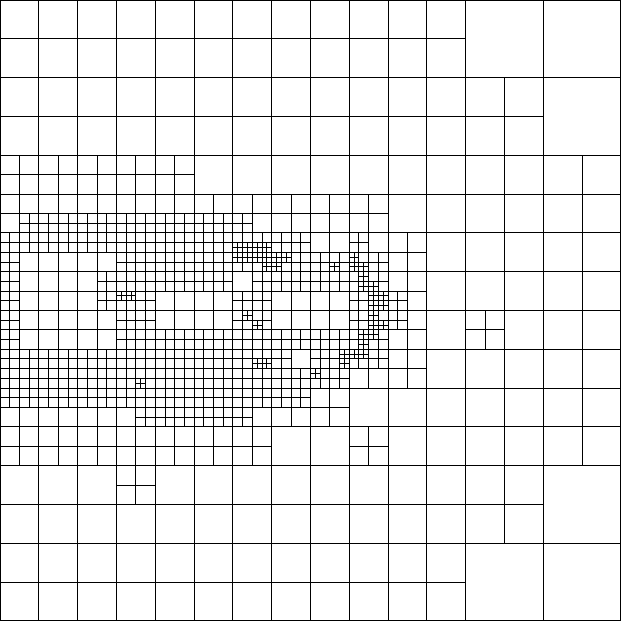}
    \hfill
    \includegraphics[width=0.16\textwidth]{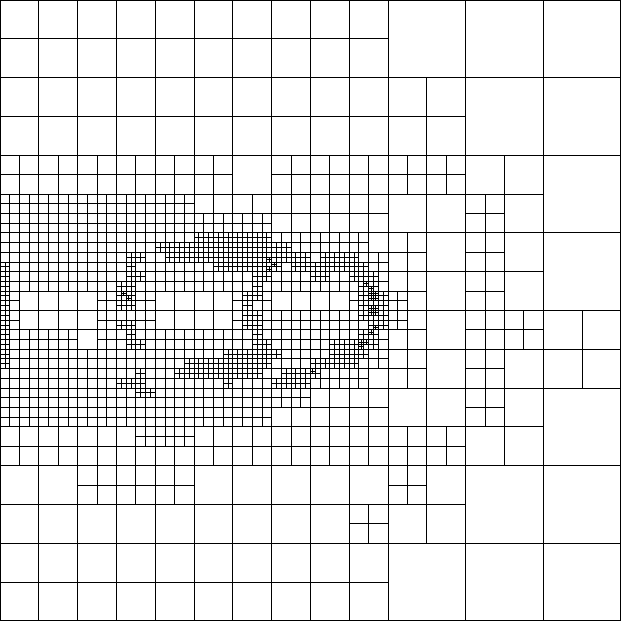}
    \caption{\it Sequences of meshes generated by different mesh refinement criteria. Top: Mesh refinement is driven by an a posteriori error indicator. Bottom: Mesh refinement is driven by a a smoothness indicator.}
    \label{fig:inverse-solution-meshes}
\end{figure}

However, we can sometimes make intuitive comparisons based on experience on ``how a good mesh should look'', even though for problems like the one under consideration, it is generally difficult to create such meshes by hand a priori. For comparison with the meshes shown in Fig.~\ref{fig:inverse-solution-information-content-mesh}, we present
in Fig.~\ref{fig:inverse-solution-meshes} the meshes generated using two different, ``traditional'' mesh refinement criteria that we will justify in more detail in Appendix~\ref{sec:amr-ee}. In both cases, the meshes are obtained by always refining those cells $K$ that have the largest ``refinement indicators'' $\eta_K$. In the top row of Fig.~\ref{fig:inverse-solution-meshes}, this indicator $\eta_K$ is the cell-wise norm of the residual of the third equation of \eqref{eq:inverse-problem-system} and is thus an a posteriori error indicator that can be derived in a way similar to that shown in \cite{Bangerth2002,BJ07ip,BR03,Griesbaum2008}:
\begin{align*}
    \eta_K = \|\beta q - \lambda\|_{L_1(K)}.
\end{align*}
We will consequently refer to this quantity as the ``error estimator''.

In the bottom row of Fig.~\ref{fig:inverse-solution-meshes}, we show meshes generated by evaluating a finite difference approximation $\nabla_h q(\mathbf x)$ on each cell by comparing the values of $q|_K$ with the values of $q$ on neighboring cells, and then computing
\begin{align*}
    \eta_K = h_K \|\nabla_h q\|_{L_2(K)},
\end{align*}
where $h_K$ denotes the diameter of cell $K$.
The choice of $\eta_K$ is proportional to the interpolation error of a continuous
function when approximated by a piecewise constant finite element
function (as we do here); the indicator therefore measures where the
piecewise constant approximation is likely poor. We will refer to this
criterion as the ``smoothness indicator''. It is also used in \cite{SambridgeRawlinson2005}, for example.

As outlined in Section~\ref{sec:vignette-mesh-refinement}, the
literature contains discussions of many other ways to refine meshes
for the inverse problem, but we consider the two mentioned above as
representative mesh refinement criteria to compare our approach against. We provide a detailed derivation of these two indicators in Appendix~\ref{sec:amr-ee}.

The meshes shown in
Fig.~\ref{fig:inverse-solution-meshes} are structurally similar to
those generated based on the information content and shown in
Fig.~\ref{fig:inverse-solution-information-content-mesh}. However,
they lack the top-bottom symmetry of the ones in
Fig.~\ref{fig:inverse-solution-information-content-mesh} and look
generally less organized, owing to the fact that \textit{they are
  based on the solution of the inverse problem, which is subject to
  the noise in the measurements}, whereas the information density
reflects only how much we \textit{know} about the solution at a
specific point in the domain -- that is, a quantity that is
\textit{independent of the concrete realization of the noise that is
  part of the measurements}, see also
Remark~\ref{rem:information-is-independent-of-noise}. Conceptually,
the best mesh should be independent of the concrete
realization of noise, although dependent on \textit{what} is being measured.  The
refinement by information content 
allows us to construct the mesh even before solving the inverse
problem because it does not depend on the solution of the inverse
problem.

If the individual measurements $m_\ell$ had had differently sized
measurement errors, then this would also have affected the information
density-based mesh refinement and led to smaller cells where more
accurate measurements are available. In contrast, there is no such
direct dependence for the other refinement methods; rather, for those methods, variable
noise levels only affect mesh refinement because $\lambda$ indirectly depends on the error level.

\subsubsection{Quantitative evaluation: Condition numbers of matrices}

A more concrete comparison between meshes would be to measure the
degree of ill-posedness of the problem. Of course, we use
regularization to make the problem well-posed, but a well-chosen mesh
results in a matrix after discretization that has a better condition
number than a poorly chosen mesh, and for which the reconstruction is
consequently less sensitive to noise. In practice, the condition
number is a poor indicator since it considers only the largest and
smallest eigenvalues; we hypothesize that a better criterion would be to ask how many
``large'' eigenvalues there are, and it is this criterion that we will
consider below.

To test this hypothesis, let us consider the discretized version of \eqref{eq:inverse-problem-system}. If we collect the degrees of freedom of a finite element discretization of $u$ into a vector $U$, and similarly those of $\lambda$ into a vector $\Lambda$ and those of $q$ into the vector $P$ (a symbol chosen to avoid confusion with the matrix $Q$ of Section~\ref{sec:finite-dimensional}), then \eqref{eq:inverse-problem-system} corresponds to the following system of linear equations after discretization by the finite element method:
\begin{align}
\label{eq:inverse-problem-system-discrete}
\begin{split}
    AU &= BP, \\
    A^T \Lambda &= -C(U - Z),\\
    \beta N P - B^T \Lambda &= 0.
    \end{split}
\end{align}
Here, the matrix $A$ corresponds to the discretized operator ${\mathcal{L}}=\mathbf b \cdot \nabla - D \Delta$ acting on the finite element space chosen to discretize the state and adjoint variables, and $N$ is the mass matrix on the finite element space chosen for the source $q$ -- that is, on the set of piecewise constant functions $s_k(\mathbf x)$ associated with the cells of the mesh. The matrix $B$ results from the product $B_{ik}=(\varphi_i,s_k)_\Omega$ between the shape functions for $u$ and $q$, and $C$ corresponds to terms of the form $C_{ij}=\sum_\ell \frac{1}{\sigma^2_\ell} \varphi_i(\mathbf \xi_\ell)\varphi_j(\mathbf \xi_\ell)$.
By noting that the matrices $A$ and $N$ are invertible, we can reduce this system of equations to an equation for $P$ by repeated substitution to
\begin{align}
  HP &= B^T A^{-1} CZ,
  \intertext{where the matrix $H$ is the Schur complement,}
  H &= B^T A^{-T} C A^{-1} B + \beta N.
\end{align}

The matrix $H$, which is symmetric and at least positive
semidefinite,  thus relates the vector of measurements $Z$ to the vector of coefficients $P$ we would like to recover. $H$ can be thought of as the discretized counterpart to the matrix $Q$ in \eqref{eq:finite-dimensional-solution}. Each eigenvalue of $H$ then corresponds to an eigenvector (``mode'') of the coefficient $q(\mathbf x)$ we would like to recover. Moreover, large eigenvalues correspond to modes that are insensitive to noise, whereas small eigenvalues correspond to modes that are strongly affected by noise. As a consequence, we would like to aim for discretizations that result in many large and few small eigenvalues.

\begin{figure}[tb]
  \includegraphics[width=0.46\textwidth]{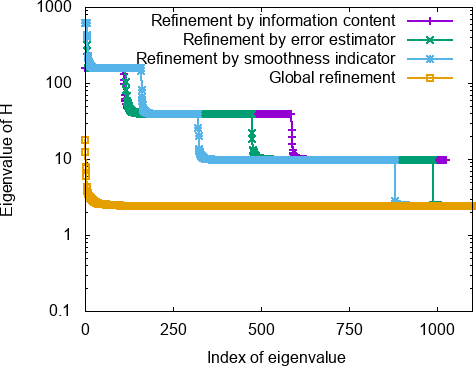}
  \hfill
  \includegraphics[width=0.46\textwidth]{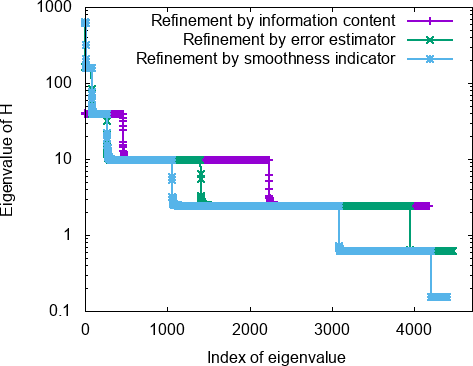}
  \caption{\it Comparison of the eigenvalues of the resolution matrix
    $H$ for different ways of refining meshes. Left: After three
    refinement cycles, yielding problems with
    approximately 1,000 parameters to identify for the adaptive refinement criteria, and 16,384 for global refinement. For the global refinement case, only the first 1,100 eigenvalues are shown; the remaining ones have nearly the same value.
    Right: After six
    refinement cycles, resulting in problems with
    approximately 4,000 parameters for the adaptive refinement criteria. (Global refinement would have resulted in more than one million parameters; the eigenvalues of this matrix could not be computed.)}
  \label{fig:eigenvalues}
\end{figure}

Fig.~\ref{fig:eigenvalues} provides a numerical evaluation of this
perspective. It shows that when refining the mesh using the
information content criterion, the eigenvalues of $H$ are further to
the top right -- in other words, there are more large eigenvalues
than when using refinement by the error estimator or the smoothness
indicator. This pattern persists after both three refinement
cycles (the left part of the figure) and six refinement cycles
(the right part).

The stair-step structure of the figure results from the fact that mesh
refinement turns one large cell into four small ones. Consequently, in
general one large eigenvalue turns into four smaller eigenvalues. By
counting the number of derivatives present in the operators that enter
into $H$, we can conjecture that the conditioning of the problem
scales with the mesh size $h$ squared; indeed, the levels in the plots
confirm that each mesh refinement step reduces the size of the
smallest eigenvalues by approximately a factor of
$(h_\text{large}/h_\text{small})^2=4$.

The left part of the figure also shows the eigenvalues of $H$ for
meshes constructed via global refinement, where the globally refined
mesh is chosen so that it has the same finest resolution as the
adaptively refined ones. Global refinement results in vastly more
unknowns than for the adaptively refined meshes, with a large majority
of small eigenvalues.  These eigenvalues all correspond to small cells
to the right of the array of detectors where very little information
is available. Given the size of the problem, we were not able to
compute the eigenvalues of $H$ following six global refinement steps;
the corresponding data are therefore omitted in the right panel.

The comparison shown in the figure confirms the hypothesis we laid
out at the beginning of the section: namely, that refining the mesh
based on information contents leads to an inverse problem that is
better posed than those that result from any of the other refinement
criteria we have compared with, in the sense that our approach leads
to more large eigenvalues.

\section{Conclusions and outlook}
\label{sec:conclusions}

In this paper, we have used a statistical approach to define how much information we have about the parameter that is recovered in an
inverse problem. More concretely, we have defined a density
$j(\mathbf x)$ that corresponds to \textit{how much we know} about the
solution at a given point $\mathbf x$, and derived an explicit
expression for it that can be computed. We have then outlined a number
of ways in which we believe that this information density can be used,
via three vignettes. Finally, we have assessed one of these
application areas numerically and showed that basing mesh refinement
for inverse problems on information densities indeed leads to meshes
that not only visually look more suited to the task than other
criteria, but also quantitatively lead to better-posed discrete
problems.

At the same time, this paper did not address many areas that would
make for very natural next questions, including the following:
\begin{itemize}
\item In our work, we have chosen a simple $L_2$ regularization term
  $\frac{\beta}{2}\|q\|_{L_2}^2$, see \eqref{eq:inverse-problem}. In
  practice, however, regularization terms would typically be used that
  penalize oscillations, for example using terms of the form
  $\frac{\beta}{2}\|\nabla^\alpha q\|_{L_2}^2$, $\alpha\ge 1$. How
  this would affect the definition of the information density would be
  interesting to ask in a future study.  In order to consider the limit
  of infinite-dimensional inverse problems, we need to ensure that the
  covariance operators that result from regularization are of trace
  class (see, for example, \cite{Stuart2010,Alexanderian2014,BuiThanh2013}), which is
  also likely a precondition for the definition of information
  densities.
    \item For many inverse problems -- such as ultrasound or seismic imaging, or electrical impedance tomography -- the quantity we would like to identify is not a right-hand source term, but a coefficient in the operator on the left side of the equation. In these cases, the definition of the information density will have to be linearized around the solution of the inverse problem, which then may make the definition of $j(\mathbf x)$ dependent on actual noise values. How this affects the usefulness of the information density is a priori unclear, but it would at least require solving the inverse problem before we can compute $j(\mathbf x)$.
      We provide a brief outlook at this case in Appendix~\ref{sec:nonlinear}.
    \item The computation of $j(\mathbf x)$ requires the solution of a
      number of forward or adjoint problems, which is expensive,
      especially for three-dimensional inverse problems, even though
      these problems are all independent of each other and can be
      computed in parallel. At the same time, although we have shown
      in Sections~\ref{sec:finite-dimensional} and
      \ref{sec:infinite-dimensional} what is necessary to form the
      complete matrix $Q$, we need only the diagonal entries of this
      matrix, see
      \eqref{eq:infinite-dimensional-information-content-3}. We may be able to compute approximations to the entries
      $Q_{kk}$ more cheaply, for example by random projections,
      low-rank approximations, or hierarchical low-rank
      approximations. Such ideas can, for example, be found in
      \cite{DasPHD,Das2010,Alexanderian2014,Spantini2015}, and in
      \cite[Section 5]{BuiThanh2013} and \cite{Ambartsumyan2020} for
      the closely related reduced Hessian matrix.
\end{itemize}
We leave the exploration of these topics to future work.

\paragraph*{Acknowledgments.}
WB would like to thank Dustin Steinhauer for an introduction to the Fisher information matrix many years ago. WB's work was partially supported by the National Science Foundation under award OAC-1835673; by award DMS-1821210; and by award EAR-1925595. CRJ's work was partially supported by the Intel Graphics and Visualization Institutes of XeLLENCE, the National Institutes of Health under award R24 GM136986, the Department of Energy under grant number DE-FE0031880, and the Utah Office of Energy Development. 

This paper describes objective technical
results and analysis. Any subjective views or opinions that might be
expressed in the paper do not necessarily represent the views of the
U.S. Department of Energy or the United States Government. Sandia
National Laboratories is a multimission laboratory managed and
operated by National Technology and Engineering Solutions of Sandia
LLC, a wholly owned subsidiary of Honeywell International, Inc., for
the U.S. Department of Energy's National Nuclear Security
Administration under contract DE-NA-0003525.


\appendix

\section{On the extension to nonlinear problems}
\label{sec:nonlinear}

The derivations in both Section~\ref{sec:finite-dimensional} (for the
finite-dimensional case) and Section~\ref{sec:infinite-dimensional} (for
the infinite-dimensional case) have assumed that the forward model is
linear, and that the solution-to-measurement map is as well. This linearity is
clearly a severe restriction -- many inverse problems are
nonlinear: typical examples are optical tomography, ultrasound
imaging, electrical impedance tomography, or seismic
imaging. In all of these cases, we seek to identify a coefficient
inside a differential operator.

In order to understand what needs to change in the formalism we have
outlined in this paper to accommodate the nonlinear case, we stay
with the finite-dimensional case for simplicity; the
infinite-dimensional case will follow the exact same steps. The
argument is essentially one of linearization, predicated
on the assumption that measurement errors are small and that the
parameter-to-measurement map is sufficiently smooth.

First, let us assume that we replace the forward model
\eqref{eq:finite-dimensional-state-eq} by the following model that is
(potentially) nonlinear in the parameters $q_k$ we seek:
\begin{align}
  A(S\mathbf q;\mathbf u) = \mathbf f,
\end{align}
where $S\mathbf q$ is a vector of coefficients that appear in the
operator and is expressed in terms of a basis encoded by the
columns of the matrix $S$ times the unknown parameters $\mathbf q$
we would like to identify. We have also allowed for a right-hand side
forcing term $\mathbf f$ that is assumed known and that could have
been incorporated into $A$ but commonly is not. 
For a given $\mathbf q$, we then use the following notation to denote
the solution of the forward model:
\begin{align}
  \mathbf u = A(S\mathbf q;\bullet)^{-1} \mathbf f.
\end{align}

Next, we assume that we obtain measurements using (possibly
nonlinear) functionals $m_\ell$ instead of the linear ones in
\eqref{eq:finite-dimensional-measurements}:
\begin{equation}
    z_\ell = m_\ell(\mathbf u) + \varepsilon_\ell.
\end{equation}

We can combine the state equation and the measurement process by
introducing the parameter-to-measurement map $\mathbf F(S\mathbf q)$:
\begin{equation}
    z_\ell = \mathbf F(S\mathbf q)_\ell + \varepsilon_\ell.
\end{equation}
For comparison with the material in
Section~\ref{sec:finite-dimensional}, note that in the linear
case, we have that
$A(S\mathbf q;\mathbf u)=A\mathbf u - S\mathbf q$, $\mathbf f=0$,
$A(S\mathbf q;\bullet)^{-1} \mathbf f = A^{-1}S\mathbf q$, and
$\mathbf F(S\mathbf q) = MA^{-1}S\mathbf q$.

With this parameter-to-measurement operator $\mathbf F$, the
deterministic parameter estimation problem
\eqref{eq:finite-dimensional-minimization-reduced} then takes on the
following form:
\begin{align}
  \label{eq:finite-dimensional-minimization-reduced-nonlinear}
  \min_{\mathbf q} \; {\mathcal{J}}_\text{red}(\mathbf q) :=&
  \underbrace{
    \frac 12 \left\| \mathbf F(S \mathbf q) - \mathbf z \right\|^2_{\Sigma^{-2}} 
    + \frac{\beta}{2} \left\| R \mathbf q \right\|^2}_{=:{\mathcal{J}}_\text{red}(\mathbf q)}.
\end{align}
Let us assume that this problem has a unique minimizer,
which we will denote by $\mathbf q^\ast=\mathbf q^\ast(\mathbf
z)$.%
\footnote{There are, of course, problems in which ${\mathcal{J}}_\text{red}$ has more than one local minimizer -- for example in
full waveform seismic wave tomography. In those cases, everything that
follows will not be useful because it relies on \textit{local}
linearization of the problem.}
Then, under appropriate smoothness assumptions, we can perform a Taylor expansion of the optimality 
conditions for $\mathbf q^\ast$. For two random noisy realizations
of measurements $\mathbf z_1,\mathbf z_2$, we obtain that the minimizers
$\mathbf q_1=\mathbf q(\mathbf z_2), \mathbf q_2=\mathbf q(\mathbf
z_2)$ of \eqref{eq:finite-dimensional-minimization-reduced-nonlinear}
have to satisfy the following relationship that generalizes
\eqref{eq:finite-dimensional-solution-difference}:
\begin{align}
    \label{eq:finite-dimensional-solution-difference-nonlinear}
    Q (\mathbf q_1 - \mathbf q_2) = S^T F'(S\mathbf q^\ast)^T \Sigma^{-2}(\mathbf z_1-
    \mathbf z_2)
    + {\mathcal{O}}(\|\mathbf z_1-\mathbf z\|^2)
    + {\mathcal{O}}(\|\mathbf z_2-\mathbf z\|^2).
\end{align}
Here, the matrix $F'(S\mathbf q^\ast)$ is the derivative of the vector function
$\mathbf F$ with regard to its argument, evaluated at $S\mathbf
q^\ast$. In the linear case,
$F'=MA^{-1}$ is independent of $q^\ast$. The matrix $Q$ above readily generalizes the situation
for the linear case and has the form:
\begin{equation}
  \label{eq:def-Q-nonlinear}
  Q=S^TF'(S\mathbf q^\ast)^{T} \Sigma^{-2} F'(S\mathbf q^\ast) S + \beta R^TR.
\end{equation}

Ultimately, the information content we seek is based on the precision
matrix, which we obtain via the Fisher information matrix from the
Bayesian inverse problem
\eqref{eq:finite-dimensional-def-prob}. Specifically, we will continue
to define the information content via
the square root of the inverse
of the diagonal elements of the covariance matrix
\eqref{eq:def-covariance-matrix},
which we approximate by way of the Cram{\'e}r-Rao bound and using \eqref{eq:fisher}:
\begin{equation}
  \label{eq:def-jk-nonlinear}
  j_k := \sqrt{(I_p)_{kk}},
  \qquad\text{ where }
  (I_p)_{kk} = {\mathbb E}\left[ \frac{\partial^2}{\partial q_k^2} {\mathcal{J}}_\text{red}(\mathbf q) \right]
\end{equation}

The differences between the linear and nonlinear cases appear
here. First, in the linear case, the Cram{\'e}r-Rao bound holds with
equality, and as a consequence the definition of $j_k$ in the previous
equation also satisfies the desired identity \eqref{eq:def-jk-via-cov},
\begin{equation*}
  j_k := \sqrt{(I_p)_{kk}} = \frac{1}{\sqrt{\text{var}_p(\mathbf q)_k}}.
\end{equation*}
However, in the nonlinear case, the probability distribution
$p(\mathbf q| \mathbf z)$ is no longer Gaussian (because ${\mathcal{J}}_\text{red}$ is no longer quadratic in $\mathbf q$), and
consequently the last identity of the previous equation may no longer
be exact.

Second, in the linear case, we had $I_p=Q$ because
${\mathcal{J}}_\text{ref}({\mathbf q})$ is a quadratic function, its second derivative is
a constant,
$\frac{\partial^2}{\partial q_k \partial q_l} {\mathcal{J}}_\text{red}(\mathbf q)=Q_{kl}$,
and consequently the expectation value in \eqref{eq:def-jk-nonlinear}
simply evaluates to $Q_{kk}$. In the nonlinear case, a
straightforward Taylor expansion of the definition of ${\mathcal{J}}_\text{red}$ around
$\mathbf q^\ast$ shows that to first order in $\mathbf q-\mathbf q^\ast$,
\begin{align*}
  I_p \approx
  Q
  +
  \text{sym}\left(
    (\mathbf F(S\mathbf q^\ast) - \mathbf z)^T\Sigma^{-2}
  F''(S\mathbf q^\ast)SS
  \right)
  +
  \text{higher order terms},
\end{align*}
with appropriate contractions over the indices of the various matrices and vectors that
appear in this expression. Here, $\text{sym}(X)=\frac{1}{2}(X+X^T)$
symmetrizes the matrix $X$.

The first correction term beyond $Q$ in the previous expression
corresponds to the difference between the Gauss-Newton and Newton
matrices in least-squares minimization \cite{NW99}. The correction is
small if either the problem is nearly linear ($F''$ is small),
or if the residual $\mathbf F(S\mathbf q^\ast) - \mathbf z$ is small
-- essentially a condition on whether the model is able to accurately
predict the data and whether the noise level is small.

The next order term of the higher order corrections in the expression above is of the form 
$
  {\mathbb E}\left[
  \text{sym}\left(
  P
  (\mathbf q-\mathbf q^\ast)
  \right)
  \right]
  =
  \text{sym}\left(
  P\;
  {\mathbb E}\left[\mathbf q-\mathbf q^\ast
  \right]
  \right)
$,
where $P$ is a rank-3 tensor composed of a lengthy list of terms proportional to at least second derivatives of $\mathbf F$, often multiplied by the residual $\mathbf F(S\mathbf q^\ast)-\mathbf z$. As above, for problems that either have a small residual or are nearly linear, $P$ by itself is already small. Furthermore, if the noise level is small, we should expect the posterior probability to be localized near $\mathbf q^\ast$, nearly symmetric around $\mathbf q^\ast$, and therefore for the expectation value ${\mathbb E}\left[\mathbf q-\mathbf q^\ast
  \right]$ to be small.

Whether or not these correction terms will be small for a concrete application depends, of course, on
both how nonlinear the operator $\mathbf F$ is (i.e., how large
$\mathbf F''$ is) and
how large we expect the measurement errors (and consequently the residuals) to be. Assuming that one or the other is indeed small for a particular case, it seems reasonable to expect that our choice
\eqref{eq:finite-dimensional-information-content} for the information
content, namely
\begin{equation*}
    j_k := \sqrt{Q_{kk}},
\end{equation*}
may still be a useful one for nonlinear problems.

Having settled the question of how to define the information content
$j_k$ in the nonlinear case, there remains the question of how to
compute it. Given \eqref{eq:def-Q-nonlinear}, we find that it can be
computed in much the same way as outlined in
Section~\ref{sec:finite-dimensional-information},
Remark~\ref{remark:compute-Q}. In analogy to the derivations there,
note that now
\begin{align*}
    Q_{kk}
    &=
    \mathbf e_k^T S^TF'(S\mathbf q^\ast)^{T} \Sigma^{-2} F'(S\mathbf
    q^\ast) S \mathbf e_k+
    \beta \mathbf r_k^T \mathbf r_k
    \\
    &=
    \left(\Sigma^{-1} F'(S\mathbf q^\ast) S \mathbf e_k\right)^T
    \left(\Sigma^{-1} F'(S\mathbf q^\ast) S \mathbf e_k\right)
    +
    \beta \mathbf r_k^T \mathbf r_k
    \\
    &=
    \sum_\ell
    \frac{1}{\sigma_\ell^2} \eta_{\ell k}^2
    +
    \beta \mathbf r_k^T \mathbf r_k,
\end{align*}
where
\begin{align*}
  \eta_{\ell k}
  = F_\ell'(S\mathbf q^\ast) S \mathbf e_k.
\end{align*}
These quantities require a solve with the linearized forward model,
linearized around the solution $\mathbf q^\ast$ of the deterministic
problem. Similar considerations apply when computing $Q_{kk}$ via
adjoint solutions (the ``alternative way'' in Remark~\ref{remark:compute-Q}).

The cost of these forward solves is the same as for the
linear case, but two other considerations come into play. First, the
solution of the nonlinear inverse problem problem (in which we are of
course primarily interested) is generally more expensive than for linear
problems, and so the cost of computing information contents $j_k$ may
or may not be a major cost any more, unlike in the linear case. On the other
hand, in the nonlinear case, the definition of $Q_{kk}$ and
consequently $j_k$ depends on the solution $\mathbf q^\ast$ of the
deterministic problem and, thus, on the noisy data $\mathbf z$ from
which $\mathbf q^\ast$ was computed. This is unfortunate: One of the
major benefits of using $j_k$ touted in
Remark~\ref{rem:information-is-independent-of-noise} was that it
did \textit{not} depend on the noisy data, unlike all other criteria we know
of. Whether this dependence is strong enough in practice to be
important will remain an open question for now; as in the discussions
about the difference between $Q$ and $I_p$, we can speculate that
this is a higher order concern, unlike for example in the mesh
refinement criteria in Fig.~\ref{fig:inverse-solution-meshes}, where
the noise in the data was clearly visible in the structure of the
generated mesh.

\section{On the derivation of traditional mesh refinement criteria}
\label{sec:amr-ee}

Section~\ref{sec:mesh-refinement-comparison} contains a comparison of our information-density based mesh refinement criterion against two other mesh refinement criteria that we have called ``error estimator'' and ``smoothness indicator''. Both of these can be derived in ways that follow common finite element theory. In the following, we will attempt to provide an intuitive derivation of their form, without attempting to be precise and complete.

The starting point for the derivation of both forms is the optimality condition \eqref{eq:inverse-problem-system} -- a system of partial differential equations -- for the deterministic problem we are trying to solve. For simplicity, let us repeat this set of equations here:
\begin{align*}
    {\mathcal{L}} u(\mathbf x) &= q(\mathbf x), \\
    {\mathcal{L}}^\ast \lambda(\mathbf x) &= -\sum_\ell \frac{1}{\sigma^2_\ell} (u(\mathbf \xi_\ell)-z_\ell) \delta(\mathbf x-\mathbf \xi_\ell), \\
    \beta q - \lambda &= 0.
\end{align*}
In order to compute a numerical solution, we seek finite-dimensional approximations $u_h(\mathbf x), \lambda_h(\mathbf x), q_h(\mathbf x)$ of the exact solution $u(\mathbf x), \lambda(\mathbf x), q(\mathbf x)$ of this problem. We will use the finite element method; specifically, we will use piecewise polynomials of degree $p_u, p_\lambda, p_q$ for this approximation. The appropriate choice in the current context is to use continuous shape functions for $u_h$ and $\lambda_h$ (because we apply a second-order differential operator, for which the weak form of the equations above requires $H^1$ solutions of which the continuous functions are a subset), and discontinuous shape functions for $q_h$ (because no differential operator is ever applied to this variable).

To assess the quality of the mesh we have chosen, we would like to quantify the ``error'', i.e., the difference between $u$ and $u_h$, and similarly for the other variables. For the second order operator $\mathcal{L}$ we consider, the appropriate norm to measure this error is
\begin{align}
\label{eq:ee-apriori}
    e := \left(
      \| \nabla (u-u_h) \|^2_{L^2(\Omega)} +
      \| \nabla (\lambda-\lambda_h) \|^2_{L^2(\Omega)} +
      \| q-q_h \|^2_{L^2(\Omega)}
    \right)^{1/2}.
\end{align}
There are now two conceptually different ways in which we can approach estimating this error: \textit{a priori} and \textit{a posteriori} error estimates. These two approaches will yield the two criteria we used in Section~\ref{sec:mesh-refinement-comparison}.

A priori error estimates are based on the premise that for many equations and under certain conditions,%
\footnote{Whether the conditions that allow us to derive such error estimates are actually satisfied for \eqref{eq:inverse-problem-system} is perhaps not actually very important here. This is because we do not seek to strictly bound the error; instead, our goal is to derive a criterion that can tell us where the error is small or large, and that can consequently be used to \textit{refine the mesh}. Experience in the finite element community has long shown that for this, less rigorous goal, \textit{heuristic} derivations of refinement criteria ``by analogy'' are sufficient and often surprisingly effective.}
one can show that on every cell $K$ of the mesh, we have
\cite{Bra97,BS02}
\begin{align*}
      \| \nabla (u-u_h) \|_{L^2(K)}
      &\le c_u \| \nabla (u-I_{h,p_u} u)\|_{L^2(K)}
      &\le C_u h_K^{p_u} \|\nabla^{p_u+1} u \|_{L^2(K)},
      \\
      \| \nabla (\lambda-\lambda_h) \|_{L^2(K)}
      &\le c_\lambda \| \nabla (\lambda-I_{h,p_\lambda} \lambda)\|_{L^2(K)}
      &\le C_\lambda h_K^{p_\lambda} \|\nabla^{p_\lambda+1} u \|_{L^2(K)},
      \\
      \| q-q_h \|_{L^2(K)}
      &\le c_q \| q-I_{h,q_q} u\|_{L^2(K)}
      &\le C_q h_K^{p_q} \|\nabla^{p_q+1} u \|_{L^2(K)},
\end{align*}
where $c_\bullet,C_\bullet<\infty$ are constants we know exist but generally do not know the value of, and $I_{h,p}$ denotes the interpolation operator from a function to a discrete, piecewise polynomial function represented on the mesh with polynomial degree $p$. The first inequality in each line above represents the ``best-approximation'' property that holds for many equations and finite element schemes, whereas the latter is the Bramble-Hilbert inequality for polynomial approximation. $h_K$ denotes the diameter of cell $K$.

Because we were primarily interested in the mesh used to discretize the coefficient $q$ in Section~\ref{sec:results-mesh}, we chose $p_u=p_\lambda=3, p_q=0$.  In other words, we used cubic elements to represent $u,\lambda$ accurately, but approximated $q$ only with piecewise constant functions. With the considerations above, the a priori error estimate \eqref{eq:ee-apriori} would then have the form
\begin{align*}
    e &= \left(
      \sum_K
      \left[
      \| \nabla (u-u_h) \|^2_{L^2(K)} +
      \| \nabla (\lambda-\lambda_h) \|^2_{L^2(K)} +
      \| q-q_h \|^2_{L^2(K)}
      \right]
    \right)^{1/2}
    \\
    &\le
     C
     \left(
      \sum_K
      \left[
      h_K^6 \|\nabla^4 u \|_{L^2(K)}^2 +
      h_K^6 \|\nabla^4 \lambda \|_{L^2(K)}^2 +
      h_K^2 \|\nabla q \|_{L^2(K)}^2
      \right]
    \right)^{1/2}
    \\
    &\approx
     C
     \left(
      \sum_K
      h_K^2 \|\nabla q \|_{L^2(K)}^2
    \right)^{1/2}.
\end{align*}
The last approximation is because for sufficiently small $h_K$, the first two terms in square brackets are much smaller than the third because their exponents are much larger. In practice, the exact solution $q$ is of course unknown, but we can derive a computable quantity for mesh refinement purposes by defining
\begin{align*}
    \eta_K = h_K \| \nabla_h q_h\|_{L_2(K)}
\end{align*}
for every cell, where $q_h$ is the computed approximate solution, and $\nabla_h q_h$ is a finite difference approximation of the gradient $\nabla q$ based on the piecewise constant function $q_h$. We have omitted the (unknown) constant $C$ because we are only interested in comparing which cells have large and which have small errors, rather than actually computing a reasonable approximation for $e$.
This derivation of $\eta_K$ therefore leads to the ``smoothness indicator'' used in Section~\ref{sec:mesh-refinement-comparison}. The term ``smoothness'' indicates that we are measuring the size of derivatives of the function $(\mathbf x)$.

The second approach toward deriving estimates for the error $e$, called ``a posteriori error estimation'', is premised on first computing a numerical solution $u_h,\lambda_h,q_h$ and then deriving an estimate for $e$ that does not require knowledge of the true solution $u,\lambda,q$. Omitting many details (but see \cite{BR03,Ver96}), we typically end up with an estimate that contains norms of the residual (i.e., the degree to which the left and right hand sides of \eqref{eq:inverse-problem-system} are not actually equal), times a constant of again unknown size, times powers of the mesh size $h_K$ in the same way as above. That is, these error estimators will have the form
\begin{align*}
    e
    &\le
     C
     \Biggl\{
      \sum_K
      \Biggl[
      h_K^6 \left(\|{\mathcal{L}}u_h-q_h \|_{L^2(K)}^2 + \text{jump}(u_h)\right) 
    \\
    & \qquad\qquad
      +
      h_K^6 \left(\left\|{\mathcal{L}}^\ast \lambda_h
         + \sum_\ell \frac{1}{\sigma^2_\ell} (u(\mathbf \xi_\ell)-z_\ell) \delta(\mathbf x-\mathbf \xi_\ell)\right\|_{L^2(K)}^2 + \text{jump}(\lambda_h)\right)
    \\
    & \qquad\qquad
      +
      h_K^2 \|\beta q_h - \lambda_h \|_{L^2(K)}^2
      \Biggr]
    \Biggr\}^{1/2}
    \\
    &\approx
     C
     \left(
      \sum_K
      h_K^2 \|\beta q_h - \lambda_h \|_{L^2(K)}^2
    \right)^{1/2}.
\end{align*}
In these formulas, $\text{jump}(u_h)$ and $\text{jump}(\lambda_h)$ are ``jump terms'' that arise due to integration by parts and whose specific form is not important in the current context. Crucially, again for reasons of high powers of $h_K$, we can estimate $e$ in the last line by omitting everything but the term that contains the residual of the last line of the optimality conditions \eqref{eq:inverse-problem-system}. Again dropping the unknown constant $C$ for the same reasons as above, we can define an indicator for the size of the error on cell $K$ by
\begin{align*}
    \eta_K = h_K \| \beta q_h - \lambda_h\|_{L_2(K)}.
\end{align*}
This is nearly the ``error estimator'' form of $\eta_K$ presented in Section~\ref{sec:mesh-refinement-comparison} -- the difference being that we use the $L_1$-norm instead of the $L_2$-norm above.%
\footnote{The choice of norm is immaterial for the derivation: At some point in the development, one has to apply H{\"o}lders inequality and can choose exponents, which can produce either of the two norms mentioned.}
We note that this derivation is essentially what all publications do that base their mesh refinement decisions on a posteriori error estimates, perhaps up to the choice of norm or using weighting factors derived from dual solutions (see, for example, \cite{Bangerth2002,BR03,Ban08ip,Griesbaum2008}, along with many others). The estimator $\eta_K$ we use here is then really just a ``redux'' of what other papers use, adapted to the fact that we have chosen $p_u,p_\lambda \gg p_q$.

The two choices for the refinement criterion $\eta_K$ we have derived in these ways thus represent the two most widely used ways to derive error estimators for finite element discretizations. There, of course, exists a large collection of publications that derive variations of such estimators, perhaps also including ways of assessing the sizes of the constants that appear in the estimates, but experience has shown that -- at least for the purposes of mesh refinement -- the meshes obtained from different estimators vary little and look conceptually very similar within each of the two broad categories of deriving these estimators. As a consequence, we consider the two choices we have presented here as representative for these categories, and we do not expect that any variation on the approach would produce meshes that are fundamentally better or worse. This, in particular, includes the methods presented in the papers mentioned in the previous paragraph. In contrast, the refinement indicator based on the information density is an entirely different approach that has no resemblance to either of the two methods outlined above -- in fact, it makes no reference to any numerical solution of the optimality conditions at all!

\bibliographystyle{siam}
\bibliography{paper}

\end{document}